\documentclass{amsart}

\usepackage{amssymb,amsmath,amscd,amsthm,xspace,color,tocvsec2, array, 
	tikz-cd}
\usepackage[mathscr]{euscript}
\usepackage{enumerate}
\usepackage[breaklinks=true]{hyperref}
\usepackage[letterpaper]{geometry}
\usepackage[all, cmtip]{xy}
\usepackage[ampersand]{easylist}
\usepackage{stmaryrd}
\usepackage[myheadings]{fullpage}
\usepackage{csquotes}
\usepackage{mathtools}

\tikzset{
	labl/.style={anchor=south, rotate=90, inner sep=.5mm}
}
\newcommand{\sS}{\mathscr{S}}

\newcommand{\DGCat}{\mathsf{DGCat}_{cont}}

\newcommand{\Fun}{\operatorname{Fun}}

\newcommand{\ch}{\operatorname{ch}}

\newcommand{\Hom}{\operatorname{Hom}}

\newcommand{\gr}{\text{gr}}

\newcommand{\g}{\widehat{\mathfrak{g}}}

\newcommand{\p}{\partial}

\newcommand{\DDGCat}{\mathsf{DGCat}}

\newcommand{\halpha}{\check{\alpha}}

\newcommand{\ft}{\mathfrak{t}}

\newcommand{\fg}{\mathfrak{g}}

\newcommand{\Fl}{\operatorname{Fl}}

\newcommand{\C}{\ensuremath{\mathbb{C}}}

\newcommand{\OO}{\mathscr{O}}

\newcommand{\Spec}{\text{Spec }}
\newcommand{\Sym}{\operatorname{Sym}}

\newcommand{\Ad}{\on{Ad}}

\newcommand{\Ext}{\operatorname{Ext}}

\newcommand{\fb}{\mathfrak{b}}
\newcommand{\gk}{\g_\kappa}
\newcommand{\fn}{\mathfrak{n}}

\newcommand{\on}{\operatorname}

\newcommand{\scc}{\mathscr{C}}
\newtheorem{theo}[subsubsection]{Theorem}
\newtheorem*{theo*}{Theorem}
\newtheorem{lem}[subsubsection]{Lemma}
\newtheorem{lemma}[subsubsection]{Lemma}
\newtheorem{cor}[subsubsection]{Corollary}
\newtheorem{pro}[subsubsection]{Proposition}

\theoremstyle{remark}
\newtheorem{defn}[subsubsection]{Definition}
\newtheorem{re}[subsubsection]{Remark}
\newtheorem{ex}[subsubsection]{Example}

\newcounter{steps}[subsection]
\newtheorem{step}[steps]{Step}

\newcommand{\normord}[1]{:\mathrel{#1}:}
\numberwithin{equation}{section}

\newcommand{\Wf}{W_{\on{f}}}

\newcommand{\fh}{\mathfrak{h}}

\newcommand{\msf}{\mathsf}

\newcommand{\sW}{\EuScript{W}}

\newcommand{\Vect}{\mathsf{Vect}}

\renewcommand{\mod}{\text{\textendash}\operatorname{mod}}
\newcommand{\Whit}{\mathsf{Whit}}

\newcommand{\xar}[1]{\xrightarrow{#1}}

\newcommand{\fc}{{\mathfrak c}}

\newcommand{\bA}{{\mathbb A}}
\newcommand{\bB}{{\mathbb B}}
\newcommand{\bV}{{\mathbb V}}
\newcommand{\bG}{{\mathbb G}}
\newcommand{\bQ}{{\mathbb Q}}
\newcommand{\bZ}{{\mathbb Z}}
\newcommand{\Ker}{\on{Ker}}
\newcommand{\Coker}{\on{Coker}}
\renewcommand{\o}[1]{\mathring{#1}}

\newcommand{\isom}{\xar{\sim}}

\newcommand{\Av}{\on{Av}}
\newcommand{\mon}{\text{\textendash}mon}
\newcommand{\heart}{\heartsuit}
\newcommand{\vac}{\on{vac}}
\newcommand{\onto}{\twoheadrightarrow}
\newcommand{\ind}{\on{ind}}
\newcommand{\Lie}{\on{Lie}}
\newcommand{\Pro}{\mathsf{Pro}}

\newcommand{\adt}{\Ad_{t^{-\check{\rho}}}}
\newcommand{\ag}{\Ad_{t^{-\check{\rho}}}\fg}

\newcommand{\sA}{\mathscr{A}}

\theoremstyle{remark}

\newcommand{\colim}{\on{colim}}
\renewcommand{\lim}{\on{lim}}

\renewcommand{\mod}{\text{\textendash}\mathsf{mod}}

\DeclarePairedDelimiter\ceil{\lceil}{\rceil}

\begin{document}

\title{Localization for affine $\sW$-algebras}
\author{Gurbir Dhillon and Sam Raskin}
\date{\today}

\maketitle
\begin{abstract}We prove a localization theorem for 
	affine $\sW$-algebras 
	in the spirit of Beilinson--Bernstein and Kashiwara--Tanisaki. 
	More precisely, for any non-critical regular weight
	$\lambda$, we identify $\lambda$-monodromic Whittaker
	$D$-modules on the enhanced affine flag variety with a full subcategory 
	of Category $\OO$ for the $\sW$-algebra.
	
	To identify the essential image of our functor, 
	we provide a new realization of Category $\OO$ for affine $\sW$-algebras
	using Iwahori--Whittaker modules for the corresponding
	Kac--Moody algebra. Using these methods, we also obtain a new proof of
	Arakawa's character formulae for simple positive energy 
	representations of the $\sW$-algebra.
\end{abstract}

\setcounter{tocdepth}{1}
\tableofcontents

\section{Introduction}

\subsection{Localization}

In their celebrated resolution \cite{beilinson-bernstein} of the 
Kazhdan--Lusztig conjecture, Beilinson and Bernstein introduced
a new method for studying irreducible modules in representation 
theory.
Their idea was to \emph{localize} representations, realizing
modules over certain rings as suitable categories of holonomic $D$-modules.
The resulting localization theorem provided Beilinson--Bernstein 
access to deep semi-simplicity results from algebraic 
geometry derived from Hodge theory and its extensions,
including Deligne's resolution \cite{weilii} 
of the Weil conjectures.

Since its emergence in \cite{beilinson-bernstein},
localization has played central role in many parts of representation 
theory, including affine
algebras and modular representations of reductive groups.

\subsection{$\sW$-algebras}

In this paper, we prove a localization theorem for affine
$\sW$-algebras. Recall that for any reductive Lie algebra 
$\fg$ equipped with a \emph{level} (i.e., 
$\Ad$-invariant symmetric bilinear form) $\kappa$,
there is an associated $\sW$-algebra denoted $\sW_{\kappa} = \sW_{\fg,\kappa}$.
These algebras have drawn interest for several reasons.

\subsubsection{} First, these algebras exhibit Feigin--Frenkel duality:
$\sW_{\fg,\kappa} \simeq \sW_{\check{\fg},\check{\kappa}}$.
Here $\check{\fg}$ is the Langlands dual Lie algebra, and
$\check{\kappa}$ is the dual level (say, for $\kappa$ 
non-degenerate). This identification is expected to 
provide the interface for quantum local geometric Langlands;
see \cite{ff-duality}, \cite{bdh}, \cite{quantum-langlands-summary}, 
and \cite{whit}.

\subsubsection{}

For $\mathfrak{sl}_2$, the algebra $\sW_{\kappa}$ is the completed enveloping 
algebra of
a Virasoro algebra with central charge depending on $\kappa$.
The Virasoro algebra is an infinite dimensional Lie algebra that
has been extensively studied due to its appearance as a symmetry
in conformal field theory and string theory. 

There is a category $\OO$ of lowest energy Virasoro representations. 
Remarkably, fairly general objects of this category appear in
physical models. For example, highest weight modules in singular
blocks appear already in free field theory, while
continuous families of Verma
modules appear in Liouville theory.  
For these reasons, lowest energy representations of the
Virasoro algebra have received considerable study.

\subsubsection{}

If $\fg$ has a simple factor of semisimple rank greater than one, 
$\sW_{\kappa}$ is no longer associated with an infinite dimensional
Lie algebra. However, there is a vertex algebra associated
with $\sW_{\kappa}$, which is conformal away from critical level.
Therefore, $\sW$-algebras provide fundamental examples
of nonlinear symmetry algebras in conformal field theory. 

\subsubsection{}

Despite considerable interest in representation theory
of $\sW$-algebras, and folklore analogies between this subject
and the geometry of affine flag variety, a direct connection
between the two has not previously been established. 
Our main result, a localization theorem for affine $\sW$-algebras,
realizes this picture in a strong sense. Even for the Virasoro algebra,
no geometric description of a category of its representations
was previously known.

\subsection{Statement of the main results}


To state our main results, we first introduce some notation 
and conventions, which are developed along with other preliminary
material in greater detail in Section \ref{s:notation}.

\subsubsection{} 

We first introduce the relevant Lie-theoretic data. 
Let $G$ be a semisimple,
simply connected group\footnote{The results we prove straightforwardly 
	extend, {\em mutatis mutandis},
	to the case of any reductive group.}  with Lie algebra $\fg$, and $\kappa$ 
	a noncritical level. Write $T$ for the abstract Cartan,  $\ft$ for its 
	Lie algebra, and $\Wf$ and $W$ for the finite and affine Weyl groups. 
Let $\lambda \in \ft^*$ be a weight that is
regular of level $\kappa$. That is, we suppose it has trivial
stabilizer under the level $\kappa$ dot action of the affine Weyl group,
cf. Section \ref{ss:regular-weights}. 

We fix once and for all a Borel subgroup $B$ with unipotent radical $N$, so 
that $B/N \simeq T$. 
In addition, we fix a nondegenerate \emph{Whittaker}
character of conductor zero for the algebraic loop group $N_F$ of $N$
\[
 \psi:N_F \to \bG_a.
\]
%

%
\subsubsection{}

The geometric side of our localization theorem is the 
DG category of $\kappa$-twisted, 
Whittaker equivariant, $\lambda$-monodromic
$D$-modules on the
enhanced affine flag variety. 
Below, we denote this category as 
\[
\Whit_{\kappa,\lambda\mon}(\Fl),
\]
cf. Section \ref{ss:monodromy} for our precise conventions. 
\subsubsection{} 
The algebraic side of the localization
theorem is the DG category of modules for $\sW_\kappa$-algebra, cf. Section 
\ref{p2:w-mods} for more detail. Below, we denote this category as
\[
\sW_{\kappa}\mod.
\]
Recall that the Verma modules $M_\chi$ for the $\sW_{\kappa}$-algebra are 
parametrized by the central characters $\chi$ of the enveloping algebra 
$\mathscr{U}(\fg)$. We identify the set of 
central characters with the set of closed points of $\Wf \backslash \ft^*$ via 
the 
Harish--Chandra homomorphism 
\[
\pi: \ft^* \rightarrow \Wf \backslash \ft^*,
\]cf. Lemma \ref{l:v2v} for our precise
normalizations. 


%

%

\subsubsection{}

With these notations in hand, we may state our main theorem. 

\begin{theo}\label{t1}

The functor of global sections on the affine flag variety yields a fully 
faithful 
embedding
\begin{equation} \label{e:t1}
\Gamma(\Fl,-):
\Whit_{\kappa,\lambda\mon}(\Fl) \to 
\sW_\kappa\mod.
\end{equation}
\noindent Moreover, its essential image is the full 
subcategory of
$\sW_{\kappa}\mod$ generated under colimits and shifts by the 
Verma modules
\[
 M_{\chi}, \quad \quad \text{for} \quad  \chi \in \pi(W \cdot 
 \lambda).
\]
That is, we obtain the subcategory generated
by Verma modules with highest weights in the image
of the affine Weyl group dot orbit of $\lambda$ under
the Harish--Chandra homomorphism. 

\end{theo}

\subsubsection{} As is typical in localization theory, the situation 
is especially nice when $\lambda$ is moreover antidominant of level $\kappa$. 
We refer to Section 
\ref{ss:regular-weights} for the standard definition of this notion, but 
emphasize that 
antidominant 
integral 
weights
exist only for negative $\kappa$.  We also recall from \cite{whit} that there 
is a canonical construction of $t$-structures on DG categories
of Whittaker equivariant $D$-modules. The construction is of
semi-infinite nature, cf. Section \ref{sss:intro-sinf}. 


\begin{theo}\label{t2}

If $\lambda$ is regular and antidominant of level $\kappa$, 
then the functor $\Gamma(\Fl,-)$ from Theorem \ref{t1} 
is $t$-exact. Moreover, this functor preserves standard,
costandard, and simple objects.
\end{theo}

In particular, passing to hearts of the
$t$-structures, we obtain a fully faithful embedding:
\begin{equation}\label{tco}
\Whit_{\kappa,\lambda\mon}(\Fl)^{\heart}
\rightarrow \OO
\end{equation}
\noindent with essential image an explicit direct sum of blocks, cf. Theorem 
\ref{t:neglev}.
Here $\OO$ is a suitable version of monodromic Category
$\OO$ for $\sW_{\kappa}$, whose block decomposition was determined
in \cite{lpw}. Indeed, a principal motivation for {\em loc. cit.} was to 
understand the essential image of our localization theorem. 

\begin{re}

As we will see in the course of the proof, both abelian
categories in \eqref{tco} are compactly generated by their subcategories of
finite length objects. Moreover, the DG category on the 
left-hand side of \eqref{e:t1} and its essential image
in the right-hand side can be canonically reconstructed from these abelian 
categories. 
Therefore, in the antidominant case, 
our results amount to an equivalence of Artinian abelian categories.

\end{re}

\begin{re}

The essential image under \eqref{tco} of the abelian category of 
$\lambda$-equivariant (not merely monodromic) $D$-modules
provides highest weight
categories of $\sW$-algebra representations at negative level.
As far as we are aware, such categories have not previously been
constructed.

\end{re}

As a consequence, we obtain the following calculation of
characters of simple modules for $\sW_{\kappa}$.



\begin{cor}\label{t3}
Let $\lambda$ be a regular antidominant weight of level $\kappa$.
Then characters of simple modules in the 
block of $\OO$ containing 
$M_{\lambda}$
are calculated via parabolic Kazhdan--Lusztig polynomials for 
\[
W_{\lambda, \on{f}} \backslash W_\lambda
\]
from the characters 
of Verma modules; here $W_\lambda$ denotes the integral Weyl group of $\lambda$ 
and $W_{\lambda, \on{f}}$ its intersection with the finite Weyl group.
\end{cor}



Corollary \ref{t3} and its extension to arbitrary blocks were
originally proved by Arakawa in the seminal paper
\cite{ara}. In the antidominant case considered here,
Corollary \ref{t3} provides a conceptual proof. Following
the original methods \cite{klic} of Kazhdan--Lusztig, 
various Kazhdan--Lusztig polynomials are well-known to 
calculate multiplicities in geometric settings. For example,
in our setting, it is essentially standard that
the Kazhdan--Lusztig polynomials considered here calculate multiplicities
in $\Whit_{\kappa,\lambda\mon}(\Fl)^{\heart}$.
Therefore, we obtain Corollary \ref{t3} as a direct consequence of 
Theorem \ref{t2}. 

In the work of Arakawa, a Drinfeld-Sokolov reduction functor from
$\gk\mod$ to $\sW_\kappa\mod$ was shown to send Vermas to Vermas
and simples to simples or zero, extending a conjecture
of Frenkel--Kac--Wakimoto \cite{fkw}. His results reduced the
calculation of simple characters for $\sW_\kappa$ 
to those for $\gk$. These
had been computed by Kashiwara--Tanisaki using localization onto
the affine flag variety,  as in the resolution of the original
Kazhdan--Lusztig conjecture \cite{ktl},\cite{kt2},\cite{kt}.
Our work provides a direct 
geometric representation theoretic explanation
for the relation between parabolic Kazhdan-Lusztig 
polynomials and $\sW$-characters.

\subsection{Relation to other work}

As far as we are aware, Theorems \ref{t1} and \ref{t2} were not explicitly
formulated as conjectures in the literature. 
However, some related results may be found. 

\subsubsection{}

For finite $\sW$-algebras, parallel results were proved by 
Ginzburg, Losev, and Webster
\cite{ginz}, \cite{los}, \cite{bw}. The expectation that
something similar should hold in affine type was mentioned in
work of Backelin--Kremnizer \cite{krem}.

\subsubsection{}

For finite $\sW$-algebras associated to a general nilpotent element $f$ of 
$\fg$,
Dodd--Kremnizer realized representations with fixed central character
in terms of asymptotic differential
operators on the standard resolution of the intersection of the
nilpotent cone with the Slodowy slice defined by $f$ \cite{dk}. 
Arakawa--Kuwabara--Malikov then gave
a chiralization of their construction to realize affine $\sW$-algebras
at critical level with fixed central character \cite{akm}. 

For a principal nilpotent $f$, as in our results, 
the intersection of the Slodowy slice and the nilpotent cone is a point. In 
addition, we work
at non-critical levels. Therefore, our 
localization results are disjoint from those of \cite{akm}.

\subsubsection{}

For $G = SL_n$, Fredrickson--Neitzke observed a
bijection between cells in a moduli spaces of wild
Higgs bundles in genus zero and minimal models for the affine 
$\sW$-algebra matching certain parameters \cite{fn}.
Our results do not speak to the phenomena they found. 

\subsection{Outline of the arguments}

Below, we first highlight the key technical aspect of the current work. We then
discuss the main features of the proofs of Theorems \ref{t1} and \ref{t2}.

\subsubsection{Semi-infinite categorical methods}\label{sss:intro-sinf}

The primary technical issue in proving Theorem \ref{t1} is 
its semi-infinite nature. As in \cite{mys}, semi-infinite 
phenomena present themselves when pairing abstract categorical methods
and infinite-dimensional Lie groups. Below, we highlight
some specific phenomena that appear in our present setting. 

Our overall strategy is to prove Theorem \ref{t1} by passing to
Whittaker equivariant objects in Kashiwara--Tanisaki's localization
theorem, which relates $D$-modules on the enhanced affine flag variety and 
Kac--Moody representations. However, there are no Whittaker equivariant 
objects in the abelian category of twisted $D$-modules on $\Fl$;
this follows because $N_F$-orbits on $\Fl$ are infinite dimensional.
Relatedly, there are no Whittaker equivariant objects in 
the abelian category of Kac--Moody representations.

However, for the corresponding \emph{DG categories} of $D$-modules and 
Kac--Moody representations, there are robust categories 
of Whittaker equivariant 
objects. By the previous paragraph, such objects are
necessarily concentrated in cohomological degree $-\infty$,
i.e., these objects are in degrees $\leqslant -n$ for all integers 
$n$.\footnote{In 
particular, as we explain in more detail below, these DG categories are {\em 
not} the unbounded derived categories of the corresponding abelian categories.}
Thus, when we form global sections of a 
Whittaker $D$-module in Theorem \ref{t1}, the 
underlying complex of vector spaces vanishes, but the
corresponding object of $\sW_{\kappa}\mod$ does not.

The connection between 
Whittaker equivariant objects in this setting and representations of
$\sW$-algebras was established by the second author in \cite{whit}.
In particular, one consequence of the main 
construction of \cite{whit}, which is crucial for the present work,
is a systematic way to realize \emph{abelian} categories of Whittaker equivariant
objects, even though such objects lie in cohomological degree $-\infty$ when 
forgetting the equivariance.

The importance of considering categories of representations
with objects in degree $-\infty$
in infinite dimensional 
algebras was first highlighted by Frenkel--Gaitsgory
\cite{fgl}. We follow them in sometimes referring to such methods
as \emph{renormalization}, and the resulting DG categories
as \emph{renormalized}. 

Therefore, to handle the described issues of semi-infinite nature, and in 
particular to pass to Whittaker equivariant objects, we use a version of 
Kashiwara--Tanisaki localization which relates the renormalized categories of 
twisted $D$-modules and Kac--Moody representations and includes equivariance 
for categorical actions of loop groups. Such an enhancement is provided by the 
first author and J. Campbell in \cite{ahc}. 

\begin{re}\label{r:sinf-loc-dream}

The dream of localization theorems  
in semi-infinite contexts is an old one, 
going back to conjectures of Lusztig and 
Feigin--Frenkel on critical level representation theory and 
localization on the semi-infinite flag manifold \cite{licm}, \cite{ffsif}.
As far as we are aware, our work is a first instance of a semi-infinite
localization theorem.

\end{re}

\subsubsection{} 

In the 
remainder of the introduction, 
we discuss some novel ingredients used to 
determine its essential image. On both sides of Theorem 
\ref{t1}, the adolescent Whittaker filtration plays an essential role.
This is a functorial filtration on Whittaker invariants 
of categorical loop group representations that was introduced 
in \cite{whit}. 

\subsubsection{}

On the geometric side, we show in Section \ref{s:awfaflgv}
that the first term of the adolescent 
Whittaker filtration suffices to describe
Whittaker $D$-modules on $\Fl$,  
i.e., we obtain a canonical equivalence between 
the Iwahori--Whittaker 
and Whittaker categories of $D$-modules on $\Fl$. 

In the untwisted setting, this result was shown in \cite{rcp}.
However, here we deduce this result from a more general
assertion. In effect,\footnote{We do not formulate the result in these
terms. The assertion stated here follows from our Theorem \ref{t:whitd}
using the methods of \cite{BZO}.} 
we show that Whittaker invariants
for a categorical representation generated by objects 
equivariant for the $n$th
congruence subgroup of $\mathfrak{L}G$ are exhausted
by the $n$th step in its adolescent Whittaker filtration, in line
with the philosophy of \cite{whit} Remark 1.22.2.

\subsubsection{}

On the representation theoretic side, the adolescent
Whittaker construction realizes subcategories of 
$\sW_{\kappa}\mod^{\heart}$ as certain categories of 
Harish-Chandra modules for $\widehat{\fg}_{\kappa}$.

We show that Category $\OO$ for $\sW_{\kappa}$ may be realized
in the first step of the adolescent Whittaker filtration, 
i.e., as Iwahori--Whittaker
modules for the Kac--Moody algebra. In Definition \ref{d:catoiw}, we
provide an explicit description of the corresponding
subcategory of Iwahori--Whittaker Kac--Moody representations; this
subcategory is mapped isomorphically 
onto Category $\OO$ for $\sW_{\kappa}$ via Drinfeld--Sokolov reduction. 

To our knowledge, the existence of such an 
Iwahori--Whittaker model for Category $\OO$ 
was not anticipated prior to \cite{whit}. 
As an application of this perspective on 
Category $\OO$, we 
give a short new proof of the character formula for simple positive energy 
representations of $\sW_{\kappa}$, originally obtained by Arakawa in the 
seminal paper \cite{ara}. 

\subsubsection{}

In order to obtain the above results, we give a new
description of the adolescent Whittaker filtration of 
$\sW_{\kappa}\mod$. 
Unlike the descriptions
in \cite{whit}, which involved Kac--Moody algebras,
our description is intrinsic to the $\sW$-algebra,
involving local nilpotency of certain explicit elements.

One important consequence of our description is that the abelian
subcategories of $\sW_{\kappa}\mod^{\heart}$ provided
by the adolescent Whittaker construction are closed
under subobjects; this was not clear to the second author
when writing \cite{whit}.

\subsubsection{} Having completed our analysis of the relevant steps in the 
adolescent Whittaker filtration, the determination of the essential image 
reduces to a problem on baby Whittaker categories. The latter question, which 
is now of finite-dimensional (i.e., no longer semi-infinite) nature, follows 
the pattern of standard 
arguments in geometric representation theory.

\subsection{Organization of the paper} In Section \ref{s:notation}, we collect 
notation
and recall preliminary material. In Section \ref{s3}, we prove a general result 
relating depth and adolescent Whittaker filtrations, and specialize it to the 
affine flag variety. In Section \ref{s4}, we give a new interpretation of the 
adolescent Whittaker filtration on $\sW_\kappa\mod$, intrisic to the action of 
$\sW_\kappa$ on representations. Building on this, in Section \ref{s5} we give 
an explicit realization of Category $\OO$ for the $\sW_\kappa$-algebra via the 
baby Whittaker category, and use it to give a new proof of the formulae for 
simple 
characters. 
In Section \ref{s6}, we obtain the localization theorem and apply the previous 
material to develop its basic properties. Finally, we include an appendix 
describing the analog of $q$-characters for our baby Whittaker model of 
Category $\OO$. This is used to establish some technical assertions invoked in 
Sections \ref{s5} and \ref{s6}.





\subsection*{Acknowledgments}

It is a pleasure to thank Roman
Bezrukavnikov, Daniel Bump, Ed Frenkel, 
Dennis Gaitsgory, Kobi Kremnizer, Ben Webster, 
David Yang, and Zhiwei Yun for their interest and encouragement.




\section{Preliminary material}\label{s:notation}

Below, we collect notation and constructions from Lie theory
and infinite dimensional geometry. The reader may wish to skip to the next 
section,
and refer back only as needed. 

\subsection{Lie theoretic notation}

\subsubsection{}

Fix a simply-connected semisimple algebraic group $G$ with 
Cartan and Borel subgroups $T \subset B$, and 
let $N$ denote the radical of $B$. Let $B^-$ be the
corresponding opposed Borel.
Let $\fg,\fb,\fn$ and $\ft$ denote the corresponding
Lie algebras. 

%

\subsubsection{The finite root system} 

We write $\Lambda$
for the weight lattice, i.e., the characters of $T$, and 
$\check{\Lambda}$ for the coweight lattice, i.e., the cocharacters of $T$. 
Within $\check{\Lambda}$, we denote the coroots by $\check{\Phi}_{\on{f}}$,
and within them the simple coroots by 
\[
\halpha_i, \quad \text{ for } i \in \mathscr{I}.
\]
%

We will write $W_{\on{f}}$ for the finite Weyl group. 
In addition to its linear action on the 
dual Cartan $\ft^\vee$, we will use its dot action. I.e., if we write $\rho$ 
for 
the half sum of positive roots in $\ft^*$,
then the dot action of an element $w$ in $W_{\on{f}}$ on $\lambda\in \ft^*$ 
is defined by
\[
w \cdot \lambda := w(\lambda + \rho) - \rho,
\]
where the right-hand action is the linear one. 

\subsubsection{Levels}

Recall that $\kappa$ denotes a \emph{level}, i.e., an 
$\operatorname{Ad}$-invariant bilinear form on $\fg$.
There are two particular levels that play distinguished roles. 
We will denote by $\kappa_c$ the \emph{critical} level, i.e., minus 
one half times the Killing form. We denote by $\kappa_b$ the 
\emph{basic} level. If $\fg$ is simple, $\kappa_b$ is 
the unique level for which the short coroots have squared
length two. I.e., for a short coroot $\check{\halpha}$, one has
\[
\kappa_b(\halpha, \halpha) = 2.
\]
For semisimple $\fg$, $\kappa_b$ is the
unique level restricting to the basic form on each simple factor.
For $G$ simple, we recall that $\kappa_c$ is necessarily a scalar multiple of
${\kappa_b}$; by definition, that scalar is minus the dual Coxeter number
of $G$.

If $\fg$ is a simple Lie algebra, we call a level $\kappa$ noncritical if 
it does not equal $\kappa_c$. We call a level $\kappa$ positive if it lies in
\[
\kappa_c + \mathbb{Q}^{\geqslant 0} \kappa_b, 
\]
and $\kappa$ negative if it is not positive. If $\fg$ is a semisimple Lie 
algebra, 
write it as a sum of simple Lie algebras
\[
\fg \simeq \bigoplus_{j \in \mathscr{J}} \fg_j 
\]
We say a level $\kappa$ for $\fg$ is noncritical if its restriction to each 
$\fg_j$ is
noncritical. Similarly, we say $\kappa$ is positive if its restriction to each 
$\fg_j$
is positive, and $\kappa$ is negative if its restriction to each $\fg_j$ is 
negative. 


%
%

\subsubsection{The affine root system} Write $F := \mathbb{C}(\!(t)\!)$ denote 
the field of Laurent series. Write $G_F$ for the algebraic loop 
group of $G$, and $\fg_F$ for its Lie algebra. Here and below, we will mean 
topological Lie algebras, so in this case 
\[
{\fg}_F := \fg \underset{\mathbb{C}} \otimes F.
\]
Associated to our level $\kappa$ is the affine Lie algebra, given as a central 
extension 
\[
0 \rightarrow \bigoplus_{j \in \mathscr{J}} \mathbb{C} \mathbf{c}_j \rightarrow 
\gk \rightarrow {\fg}_F \rightarrow 0. 
\]
Explicitly, for elements $X$ and $Y$ of $\fg_i$ and $\fg_j$, for $i$ and $j$ in 
$\mathscr{J}$, and Laurent series $f$ and $g$ in $\C(\!(t)\!)$, the bracket
is given by
\[
[X \otimes f, Y \otimes g ] \hspace{1mm} = \hspace{1mm} [X,Y] \otimes 
fg\hspace{1mm} + \hspace{1mm} \delta_{i,j} \cdot \kappa(X, Y) 
\cdot {\on{Res}} \hspace{1mm} df g \cdot \mathbf{c}_i,
\]
where $\on{Res}$ denotes the residue, and $\delta_{i,j}$ the Kronecker delta 
function. 

%
Consider the affine Cartan
\[\ft \oplus \bigoplus_{j \in \mathscr{J}} \mathbb{C} \mathbf{c}_j.
\]
We may write its linear dual as 
\begin{equation} \label{e:dualaffcartan}
\ft^* \oplus \bigoplus_{j \in \mathscr{J}} \mathbb{C} \mathbf{c}_j^*,
\end{equation}
where $\mathbf{c}^*_j$ pairs with
$\ft$ and the $\mathbf{c}_i$, for $i \neq j$, by zero and pairs with 
$\mathbf{c}_j$ by one. 
Let us denote the real affine coroots by $\check{\Phi}$, and the simple affine 
coroots by 
\[
\halpha_i, \quad \text{ for } i \in \hat{\mathscr{I}}.
\]
Explicitly, $\widehat{\mathscr{I}} := \mathscr{I} \sqcup \mathscr{J}$. To write 
the corresponding
elements of $\ft \oplus \mathbb{C} \mathbf{c}$, if we denote by 
$\check{\theta}_j$ the short
dominant coroot of $\fg_j$, by $\kappa_j$ the restriction of $\kappa$ to 
$\fg_j$, and $\kappa_{b, j}$
the basic level for $\fg_j$, then they are given by
\[
\halpha_i, \quad \text{for } i \in \mathscr{I}, \quad \quad \text{ and } \quad 
\quad \halpha_j := -\check{\theta}_j +  
\frac{\kappa_j}{\kappa_{b,j}} \hspace{1mm}\mathbf{c}_j, \quad \text{for } j \in 
\mathscr{J}.
\]
We will write $W$ for the affine Weyl group. In addition to its linear action 
on the dual affine Cartan, we will use its dot action. I.e., if we write 
$\check{\rho}$ for the unique element of \eqref{e:dualaffcartan} satisfying
\[
\langle \check{\rho}, \halpha_i \rangle = 1, \quad \text{for } i \in 
\widehat{\mathscr{I}},
\]
then the dot action of $w$ in $W$ on an element $\lambda$ of 
\eqref{e:dualaffcartan} is given by
\[
w \cdot \lambda = w( \lambda + \hat{\rho}) - \hat{\rho}.
\]
We may identify $\ft^*$ with the affine subspace 
\[
\ft^* + \sum_{j \in \mathscr{J}} \mathbf{c}_j^*. 
\]
This affine subspace is preserved by both the linear and dot actions of $W$, 
and we will 
only need these induced affine linear actions on $\ft^*$ in what follows.

\subsubsection{Weights and integral Weyl groups} \label{ss:regular-weights} 
Recall that a weight $\lambda$ of $\ft^*$ is 
	\emph{antidominant} if
\[
\langle \lambda + \check{\rho}, \halpha_i \rangle \notin \mathbb{Z}^{> 0}, 
\quad \text{for } i \in \widehat{\mathscr{I}},
\]
and that $\lambda$ is \emph{dominant} if 
\[
\langle \lambda + \check{\rho}, \halpha_i \rangle \notin \mathbb{Z}^{< 0}, 
\quad \text{for } i \in \widehat{\mathscr{I}}.
\]
If $\kappa$ is positive, then the dot orbit of $\lambda$ always contains a 
dominant weight,
and if $\kappa$ is negative, then the dot orbit of $\lambda$ always contains an 
antidominant weight.

Associated to any weight $\lambda$ of $\ft^*$ is the subset of integral real 
affine coroots
\[
\check{\Phi}_\lambda := \{ \halpha \in \check{\Phi}: \langle \halpha, \lambda 
\rangle \in \mathbb{Z} \}.
\]
The integral Weyl group $W_\lambda$ is the subgroup of $W$ generated by the 
reflections 
indexed by the integral coroots.

We recall that a weight $\lambda$ is said to be regular
if its stabilizer in $W$ under the dot action is trivial. If $\kappa$ is not 
critical, this coincides with its stabilizer within $W_\lambda$ under the dot 
action.




%





\subsection{The $\sW_{\kappa}$-algebra} 

\subsubsection{} For $\gk$ the affine algebra as above,
we let $\sW_\kappa$ denote the $\sW$-algebra associated to $\gk$
and a principal nilpotent element $f$ of $\fg$. I.e., $\sW_{\kappa}$ is the vertex 
algebra obtained as
the quantum Drinfeld-Sokolov reduction of the vacuum representation
of $\gk$, cf. \cite{aran} and Chapter 15 of \cite{fbz}. For $\kappa$ 
noncritical, $\sW_{\kappa}$ contains a canonical conformal vector. 
Throughout this paper, unless otherwise specified, we assume $\kappa$ is 
noncritical. 

\subsubsection{}\label{p2-zhu} We recall that the Zhu algebra of $\sW_{\kappa}$ 
identifies with the center of the universal enveloping algebra
of $\fg$. We refer to Lemma \ref{l:v2v} for our
precise normalization of this isomorphism. In 
particular, for a central character 
$\chi$, we denote by
$M_\chi$ the corresponding Verma module for $\sW_{\kappa}$, and
by ${L}_\chi$ its unique simple quotient.



\subsection{Categories of representations}

\subsubsection{}
For a DG category $\scc$ equipped with a $t$-structure, we write 
$\scc^{\heartsuit}$ for the 
abelian category of objects lying in its heart, and $\scc^+$ for its full 
subcategory consisting
of bounded below, i.e., eventually coconnective, objects.

\subsubsection{} We denote the abelian category of smooth modules for $\gk$ on 
which each
$\mathbf{c}_j$, for $j \in \mathscr{J}$, acts via the identity by
\begin{equation} \label{e:gkheart}
\gk\mod^\heartsuit.
\end{equation}
We denote by $\gk\mod$ the DG category of smooth modules for $\gk$, 
which is a renormalization of the unbounded DG derived category of
$\gk\mod^\heartsuit$ introduced by Frenkel--Gaitsgory \cite{fg}. 
Explicitly, within the bounded derived category of \eqref{e:gkheart}, 
consider the subcategory 
generated under cones and shifts from modules
induced from trivial representations 
of compact open subalgebras of $\gk$;
by definition, $\gk\mod$ is the ind-completion of this 
category. The category $\gk\mod$ carries a $t$-structure 
with heart \eqref{e:gkheart}; its
bounded below part identifies canonically 
with the bounded below derived category
 of \eqref{e:gkheart}. Further details may be found in Sections 22 and 23
 of \cite{fg} as well as the very readable Section 2 of \cite{gn}. 

\subsubsection{} \label{p2:w-mods} We denote the abelian category of modules, 
in the sense
of vertex algebras, over $\sW_{\kappa}$ by
\begin{equation}
\sW_{\kappa}\mod^\heartsuit. \label{e:wmodh}
\end{equation}
We denote by $\sW_{\kappa}\mod$ the DG category of modules for $\sW_{\kappa}$.
This is a renormalization of the unbounded DG derived category of 
$\sW_{\kappa}\mod^\heartsuit$
introduced in \cite{whit}. As in the Kac--Moody case, it 
may be constructed as the ind-completion of the subcategory
generated by an explicit 
collection of modules within the bounded derived category of \eqref{e:wmodh}, 
cf. 
Section 4 of {\em loc. cit. } and Section \ref{s4} of the present paper. 

\subsection{Conventions on DG categories}

\subsubsection{}

We denote by $\DDGCat$
the $(\infty, 1)$-category of DG categories and DG functors
between them. We let $\DGCat \subset \DDGCat$ 
denote the $1$-full subcategory of 
of cocomplete DG categories and continuous DG functors between
them. We refer to \cite{gaitsroz} for background material
on the $\infty$-categorical perspective on DG categories.



\subsection{D-modules}

\subsubsection{} For a general treatment of $D$-modules on 
infinite dimensional varieties, we refer the reader to 
\cite{rdm} and \cite{beraldo}. However, with the 
exception of 
($\kappa$-twisted)
$D$-modules on the loop group, we will only deal with the DG
category of $D$-modules on ind-schemes $X$ of ind-finite type. 

Let us summarize the relevant aspects of the theory in this
simpler case. For an indscheme $X$ written\footnote{Of course,
for cardinality reasons, 
not every ind-scheme can be written as such a union. 
However, this is satisfied for our examples of interest, so for
concreteness, we assume it.}
as an ascending union 
\[
Z_0 \rightarrow Z_1 \rightarrow Z_2 \rightarrow \cdots 
\]
of schemes of finite type along closed embeddings, the category
$D(X)$ of $D$-modules on $X$ is the colimit in $\DGCat$
of the corresponding categories for the $Z_i$ along $*$-pushforwards, i.e., 
\begin{equation}
D(X) \simeq \colim D(Z_0) \rightarrow D(Z_1) \rightarrow D(Z_2) \rightarrow 
\cdots .
\end{equation}
In particular, by the $t$-exactness of these pushforwards, 
$D(X)$ has a natural $t$-structure. 
Moreover, its bounded below part $D(X)^+$ 
identifies
with the bounded below derived category of its heart $D(X)^\heartsuit$, cf. 
Lemma
5.4.3 of \cite{whit}. Finally, we should note that $D(X)$ is compactly 
generated, with
\[
D(X)^c \simeq \colim D(Z_0)^c \rightarrow D(Z_1)^c \rightarrow D(Z_2)^c 
\rightarrow \cdots,
\]
where the superscript $c$ denotes compact objects and the appearing colimit is 
taken in 
$\DDGCat$. Plainly, 
a compact object of $D(X)$ is simply a bounded complex of $D$-modules
with coherent cohomology $*$-extended from some $Z_n$, and $\Hom$ between two 
such objects may 
be computed in any $Z_n$ containing both their supports. 


\subsubsection{}

We recall that $\kappa$ defines 
categories of twisted $D$-modules the loop group
$G_F$ and its descendants. More precisely, there is a monoidal DG category
$D_{\kappa}(G_F)$ defined in \cite{whit} Section 1.30
and \cite{mys}; in the latter source,
this category is denoted $D_{\kappa}^*(G(K))$. 
This twisting is canonically 
trivialized on the 
compact open subgroup $G_O$ of regular maps from the
disc to $G$, and below we use the restriction of this trivialization to $I 
\subset G_O$. Similarly, the twisting is canonically trivialized on the group 
ind-scheme
$N_F \subset G_F$ of loops into $N$.

\subsection{Group actions on categories}

\subsubsection{} Throughout this paper, we make essential use of techniques 
from categorical representation theory. We review below some basic features we 
employ; the reader may wish to consult the user-friendly \cite{beraldo} or the 
foundational paper \cite{mys} as further references.
We also refer to \cite{paris-notes} for some useful lecture
notes on the subject.

\subsubsection{} Given a group ind-scheme $H$, by functoriality the 
multiplication on $H$ endows the category $D(H)$ with a monoidal structure 
given by convolution. We denote the associated $(\infty-)$category of DG 
categories equipped with an action of $D(H)$ by 
\[
D(H)\mod := D(H)\mod(\DGCat). 
\]
\subsubsection{} Given a $D(H)$-module $\scc$, one may form its categories of 
invariants and coinvariants, respectively given by
\[
\scc^H := \msf{Hom}_{D(H)\mod}(\msf{Vect}, \scc) \quad \text{and} \quad \scc_H 
:= \msf{Vect} \underset{D(H)} \otimes \scc. 
\]

\begin{ex} Suppose for simplicity that $H$ is of finite type, and acts on a 
scheme $X$. This endows $D(X)$ with an action of $D(H)$ by convolution. Using 
the bar resolution, one may identify $D(X)^H$ with the category of 
$H$-equivariant $D$-modules on $X$, i.e. the category of $D$-modules on the 
stack $X/H$. For a more general discussion including infinite type see 
\cite{rdm}, particularly Section 6.7. 
\end{ex}

There are tautological `forgetting and averaging' adjunctions 
\[
\on{Oblv}: \scc^H \rightleftarrows \scc: \on{Av}_{H, *}
\]
\[
\on{ins}^L: \scc_H \rightleftarrows \scc: \on{ins}.
\]

\subsubsection{} If $H$ is a group scheme whose pro-unipotent radical is of 
finite codimension, then there is a canonical identification of invariants and 
coinvariants. Namely, the map $\scc_H \rightarrow \scc^H$ induced by 
$\on{Av}_{H, *}$ is an equivalence. We will use these equivalences implicitly 
throughout.

\subsubsection{} We will also need the twists of the above by a character. The 
data of a multiplicative $D$-module $\chi$ on $H$,\footnote{Strictly speaking, 
in this generality this should be regarded as an object of the {\em dual} 
category $D^!(H)$.}, is equivalent to an action of $D(H)$ on $\msf{Vect}$ 
which we denote by $\msf{Vect}_\chi$. One has associated categories of twisted 
invariants and coinvariants 
\[
\scc^{H, \chi} := \msf{Hom}_{D(H)\mod}(\msf{Vect}_\chi, \scc) \quad \text{and} 
\quad \scc_{H, \chi} := \msf{Vect}_\chi \underset{D(H)} \otimes \scc. 
\]
These fit into adjunctions as above, and under the same hypotheses on $H$ one 
may canonically identify twisted invariants and coinvariants.

\subsubsection{} We will need two examples of nontrivial multiplicative 
$D$-modules. First, given an element $$\lambda \in \Hom(H, \mathbb{G}_m) 
\underset{\mathbb{Z}} \otimes \mathbb{C}$$
one has an associated multiplicative $D$-module $``t^{\lambda}"$ on $H$. We 
denote the corresponding twisted invariants or coinvariants by a superscript or 
subscript $(H, \lambda)$, respectively. In what follows, 
this will principally be applied to Iwahori subgroups of $G_F$. 

Second, given an element $\psi \in \on{Hom}(H, \mathbb{G}_a)$ one has an 
associated multiplicative $D$-module $``e^\psi"$ on $H$. We denote the 
corresponding twisted invariants and coinvariants by a superscript or subscript 
of $(H, \psi)$, respectively. This will principally be applied 
to $N_F$ and related prounipotent subgroups of $G_F$.

\subsection{Whittaker models and the adolescent Whittaker filtration}

\subsubsection{} Applying the general constructions of the previous section, to 
a $D_\kappa(G_F)$-module $\scc$ one may attach its Whittaker invariants and 
coinvariants 
\[
  \scc^{N_F, \psi} \quad \text{and} \quad \scc_{N_F, \psi}.
\]
A principal result of \cite{whit} is that these may be canonically identified. 
This is non-trivial because $N_F$ is not a group scheme, 
but rather a group 
ind-scheme. We employ this identification throughout. 
Unless an argument is biased toward one of these 
perspectives, we refer to them both as the 
{\em Whittaker model} of 
$\scc$, and we denote this category by $\msf{Whit}(\scc)$. 

\subsubsection{} 
Refining the equivalence between invariants and
coinvariants above, 
\cite{whit} constructed a canonical filtration of any 
Whittaker model by full subcategories
\[
   \msf{Whit}^{\leqslant 1}(\scc) \subset \msf{Whit}^{\leqslant 2}(\scc) 
   \subset 
   \cdots \quad \quad \underset{n}{\colim} \, \msf{Whit}^{\leqslant n}(\scc) \simeq 
   \msf{Whit}(\scc).
 \]
As we presently review, these may be identified with the twisted invariants for 
certain compact open subgroups $\mathring{I}_n$ of $G_F$.

\subsubsection{} Fix a positive integer $n$. Consider the $n^{th}$ congruence 
subgroup $K_n$ of $G_O$, and write $G_{n-1}$ for the quotient, and similarly 
consider $N_{n-1}$. One may form the subgroup 
of $G_O$ consisting of arcs which until $n$th order lie in $N$, i.e. 
\[ \mathring{J}_n := G_O \underset{G_{n-1}} \times  N_{n-1}.
\]
One then obtains $\mathring{I}_n$ by conjugation, namely 
\[
\mathring{I}_n := \on{Ad}_{t^{-n \check{\rho}}} \mathring{J}_n.
\]
\subsubsection{}
Each $\mathring{I}_n$ admits a unique 
additive character $\psi$ which is (i) trivial on $\mathring{I}_n \cap B^-_F$ 
and (ii) agrees with the Whittaker character on $\mathring{I}_n \cap N_F$. 
There are canonical identifications 
\[
 \msf{Whit}^{\leqslant n}(\scc) \simeq \scc^{\mathring{I}_n, \psi},
\]
and the transition functors in the above filtration are given by averaging, cf. 
Section 2 of \cite{whit} for more details. 

\subsubsection{} The first case step in the filtration will be of particular 
importance 
to us. Note that $\mathring{I}_1$ is the prounipotent 
radical of an Iwahori subgroup. For ease of notation, we often denote 
them by $\mathring{I}$ and $I$, respectively. In what follows, we sometimes 
refer to the corresponding invariants as the {\em baby Whittaker model}.

\section{Depth and adolescent Whittaker}
\label{s3}
\label{s:awfaflgv}
The goal of this section is to show that the Whittaker category on the enhanced 
affine flag variety is exhausted by the first step in its adolescent Whittaker 
filtration, i.e. canonically identifies with the baby Whittaker category. 

 We will deduce the above statement from the following more 
 general assertion. For a nonnegative integer $n$, denote by $K_n$ the 
 corresponding congruence 
 subgroup of $G_O$.

\begin{theo}\label{t:whitd} The adolescent Whittaker filtration for 
$D_\kappa(G_F/K_n)$ is exhausted by its $n^{th}$ step, i.e. 
	\[
	 \iota_n: \mathsf{Whit}^{\leqslant n} (D_\kappa(G_F / K_n)) \simeq 
	 \mathsf{Whit}( 
	 D_\kappa(G_F / K_n)). 
	\]
\end{theo}

\begin{proof} As $\iota_n$ is fully faithful, it suffices to show its essential 
surjectivity. Recalling that $\check{\Lambda}$ denotes the cocharacter lattice 
of $T$, the Iwasawa decomposition yields a stratification of 
$G_F/K_n$ with strata
	\[
	 \mathscr{S}^{\check{\lambda}} := N_F \hspace{.5mm} 
	 t^{\check{\lambda}}\hspace{.5mm} G_O/K_n, \quad 
	 \quad \text{for} \quad \lambda \in 
	 \check{\Lambda}.
	\]
Using the normality of $K_n$ in $G_O$, a standard argument shows that the 
Whittaker category on the stratum vanishes, i.e. that 
\[
\msf{Whit}(D_\kappa(\mathscr{S}^{\check{\lambda}})) \simeq 0,
\]
unless $\check{\lambda} + 
n \check{\rho}$ is a dominant cocharacter. It follows that the closure of each 
stratum contains only finitely many other strata which support Whittaker 
sheaves. In particular, to see the essential surjectivity of $\iota_n$ it 
suffices to consider objects $*$-extended from a single stratum. 
	
	To prove the latter claim, fix a cocharacter $\check{\lambda}$ such that 
	$\check{\lambda} + n \check{\rho}$ is dominant. Under this assumption, 
	using the triangular decomposition of $\mathring{I}_n$ it is 
	straightforward to see that 
	\begin{equation}
	 \mathring{I}_n\hspace{.5mm} t^{\check{\lambda}}\hspace{.5mm} G_O / K_n = 
	 (\mathring{I}_n \cap N_F) 
	 \hspace{.5mm} 
	 t^{\check{\lambda}} \hspace{.5mm} G_O/K_n,  \label{e:ids}
	\end{equation}
	i.e. that these coincide as locally closed sub-ind-varieties of $G/K_n$. 
	Similarly, our assumption on $\check{\lambda}$ implies that 
	\begin{equation} \label{e:stabs}
	 (N_F \cap t^{\check{\lambda}} G_O t^{-\check{\lambda}}) \hspace{.5mm} 
	 \subset \hspace{.5mm}  
	 (\mathring{I}_n \cap N_F).  
	\end{equation}
	It follows from \eqref{e:stabs} that we have
	\begin{equation}
	 \mathscr{S}^{\check{\lambda}} \simeq N_F \overset{N_F \cap \mathring{I}_n} 
	 \times \mathring{I}_n \hspace{.5mm} t^{\check{\lambda}} 
	 G_O/K_n.\label{e:iws2}
	\end{equation} 
	Recalling that $\mathring{I}_n$ is prounipotent, we deduce equivalences
	\begin{equation}
	 D_\kappa(\mathring{I}_n \hspace{.5mm} t^{\check{\lambda}} \hspace{.5mm} 
	 G_O/K_n)_{\mathring{I}_n, \psi} \overset{\eqref{e:ids}}\simeq 
	 D_\kappa(\mathring{I}_n 
	 \hspace{.5mm} t^{\check{\lambda}} \hspace{.5mm} G_O/K_n)_{\mathring{I}_n 
	 \cap N_F, \psi} \overset{\eqref{e:iws2}}\simeq D_\kappa( 
	 \mathscr{S}^{\check{\lambda}})_{N_F, 
	 \psi}.   \label{e:iweq}
	\end{equation} 
	It is straightforward to see that $*$-pushforward to $G_F/K_N$ intertwines 
	\eqref{e:iweq} and $\iota_n$, which yields the claimed essential 
	surjectivity. \end{proof}

We now deduce the desired consequence for the affine flag variety. 

\begin{cor} \label{c:gendeg1}The adolescent Whittaker filtration on 
$\msf{Whit}_{\kappa, 
\lambda-mon}(\Fl)$ is exhausted by its first step, i.e. 
	\[
	 \iota_1: \msf{Whit}^{\leqslant 1}_{\kappa, \lambda-mon}(\Fl) \simeq 
	 \msf{Whit}_{\kappa, \lambda-mon}(\Fl).
	\]
\end{cor}

\begin{proof} We first note that a fully faithful embedding of 
$D_\kappa(G_F)$-modules $\mathring{\scc} \rightarrow \scc$ induces equivalences
	\begin{equation} \label{e:awints}
	      \msf{Whit}^{\leqslant n}(\mathring{\scc}) \simeq ( \mathring{\scc} 
	      \cap \msf{Whit}^{\leqslant n}(\scc)),  
	\end{equation}
	i.e., they coincide as full subcategories of $\scc$. Indeed, this claim 
	follows tautologically from the identification, for any 
	$D_\kappa(G_F)$-module $\sS$, 
	\[
	  \msf{Whit}^{\leqslant n}(\sS) \simeq \msf{Hom}_{D_\kappa(G_F)\mod}( 
	  D_\kappa(G_F)^{\mathring{I}_n, \psi}, \sS ).
	\]
	To use this, note that pullback yields a fully 
	faithful embedding
	\begin{equation}
	\label{e:pullback1stK}
	   D_\kappa(G_F)^{I, \lambda-mon} \rightarrow D_\kappa(G_F)^{K_1}. 
	\end{equation}
	By Theorem \ref{t:whitd}, the adolescent Whittaker filtration for the 
	right-hand side of \eqref{e:pullback1stK} is exhausted by its first step, 
	hence by \eqref{e:awints} the same holds for the left-hand side, as 
	desired. 
\end{proof}

\section{The adolescent Whittaker filtration on representations of the 
$\sW_\kappa$-algebra revisited}
\label{s4}
\label{s:adol}

The goal for this section is to show the following result.

\begin{theo}\label{t:whit-n-sub}
	
	The subcategory $\Whit^{\leqslant n}(\widehat{\fg}_{\kappa}\mod)^{\heart}
	\subset \sW_{\kappa}\mod^{\heart}$ 
	is closed under subobjects.
	
\end{theo}


This is {\em a priori} nonobvious because 
the definition of adolescent Whittaker 
filtration is in terms of Harish--Chandra modules for varying compact open 
subgroups of $G_F$ which are related by averaging functors. We will prove this 
result by realizing the subcategory
$\Whit^{\leqslant n}(\widehat{\fg}_{\kappa}\mod)^{\heart}$
in more explicit terms, intrinsic to the action of 
the $\sW_\kappa$-algebra on a 
module. 

\subsection{Background}

Recall that \cite{whit} 
introduced the generalized vacuum modules, i.e. a 
sequence of representations 
$\sW_{\kappa}^n \in \sW_{\kappa}\mod^{\heart}$. We presently review their basic 
properties. 

We remind that these
modules have distinguished \emph{vacuum}
vectors $\vac_n \in \sW_{\kappa}^n$,
and there are natural surjections
$$\alpha_n:\sW_{\kappa}^{n+1} \onto \sW_{\kappa}^n$$
which send $\vac_{n+1}$ to $\vac_n$. More generally, 
for $m \geqslant n$, we have a composition which we denote by
\begin{equation}\alpha_{n,m}:\sW_{\kappa}^m \onto 
\sW_{\kappa}^n.\label{e:alphas}\end{equation}
As constructed in \cite{whit}, each module $\sW_{\kappa}^n$ carries a canonical
\emph{Kazhdan-Kostant} filtration. We denote these
filtrations by $F_{\bullet}^{KK} \sW_{\kappa}^n$. These satisfy $$F_{-1}^{KK} 
\sW_{\kappa}^n = 0 \quad \text{and} \quad F_0^{KK} \sW_{\kappa}^n = k \cdot 
\vac_n.$$The morphisms \eqref{e:alphas} are strictly compatible with the 
Kazhdan--Kostant filtrations.


Finally, we recall that the underlying vector space of each $\sW_{\kappa}^n$ 
has a grading
defined by loop rotation. We denote the $j$th graded
piece of $\sW_{\kappa}^n$ by $\sW_{\kappa}^n(j)$.
More importantly, for any $j$, define the subspace
$$\sW_{\kappa}^n(\geqslant j) \coloneqq \bigoplus_{k \geqslant j}
\sW_{\kappa}^n(k).$$

\subsection{Main construction}

Suppose $m \geqslant n$ are nonnegative integers. We will
construct generators for $\ker(\alpha_{n,m})$
in terms of the above data. To do so, define
\[
V_{n,m} \coloneqq
\sum_{i>0} \big(F_i^{KK} \sW_{\kappa}^{m} \cap \sW_{\kappa}^{m}
(\geqslant n\cdot i)\big) \subset \sW_{\kappa}^m.
\]
The following result plays a key role.

\begin{lem}\label{l:alpha-n-m}
	
	$V_{n,m}$ is contained in $\ker(\alpha_{n,m})$.
	Moreover, $V_{n,m}$ generates this kernel,
	i.e., $\ker(\alpha_{n,m})$ is the minimal
	subobject of $\sW_{\kappa}^m$ in $\sW_{\kappa}\mod^{\heart}$
	containing $V_{n,m}$.
	
\end{lem}

\begin{proof}We break the proof into several steps. 
	
	\label{st:alpha-3}

	\step First, we explain that the analogous
	result for Kac--Moody algebras is straightforward. Let $\fh$ be a reductive 
	Lie algebra,
	and let $\kappa_{\fh}$ be a level for $\fh$.
	
	The analogues of the generalized vacuum modules are given by 
	$$\bV_{\fh,\kappa_{\fh}}^n \coloneqq
	\ind_{t^n\fh_O}^{\widehat{\fh}_{\kappa_{\fh}}}(k).$$
	These module are equipped with PBW filtrations
	and gradings defined by loop rotation. We have
	evident structure morphisms
	$$\beta_{n,m}:\bV_{\fh,\kappa_{\fh}}^m \to
	\bV_{\fh,\kappa_{\fh}}^n$$for $m \geqslant n$ that are
	strictly compatible with the PBW filtration and compatible
	with the grading.
	
	Define $W_{n,m} \subset \bV_{\fh,\kappa_{\fh}}^m$
	by analogy with $V_{n,m}$. We claim that it is
	contained in $\ker(\beta_{n,m})$ and generates it
	as a Kac--Moody representation. Indeed, to show the containment
	it suffices to show for every positive integer $i$ that 
	$$(F_i^{PBW} \bV_{\fh,\kappa_{\fh}}^m \cap \bV_{\fh,\kappa_{\fh}}^m
	(\geqslant n\cdot i)) \hspace{.5mm} \subset \hspace{.5mm} 
	\ker(\beta_{n,m}).$$ 
	We will proceed by induction on $i$. For $i = 1$, we have a natural 
	identification
	$$F_1^{PBW} \bV_{\fh,\kappa_{\fh}}^m 
	\simeq
	\widehat{\fh}_{\kappa_{\fh}}/t^m\fh_O.$$Clearly, the
	subspace of vectors of degree $\geqslant n$ in the left-hand side 
	identifies with	$t^n\fh_O/t^m\fh_O$, and hence lies in
	$\ker(\beta_{n,m})$ by definition.
	
	Now inductively, fix $i > 1$ and consider an element  $$\xi \in
	F_i^{PBW} \bV_{\fh,\kappa_{\fh}}^m \cap \bV_{\fh,\kappa_{\fh}}^m
	(\geqslant n\cdot i).$$By definition, its symbol $$\sigma(\xi) \in 
	\gr_i^{PBW} \bV_{\fh,\kappa_{\fh}}^m =
	\Sym^i(\fh_F/t^m\fh_O)$$has degree $\geqslant n \cdot i$,
	which readily implies that it lies in
	$\Sym^{i-1}\big(\fh_F/t^m\fh_O\big) \cdot
	t^n\fh_F/t^m\fh_O.$ Therefore, there
	exists $\widetilde{\xi} \in
	F_{i-1}^{PBW} U(\widehat{\fh}_{\kappa_{\fh}}) \cdot t^n\fh_O$
	with $\xi - \widetilde{\xi} \cdot \vac_n \in
	F_{i-1}^{PBW} \bV_{\fh,\kappa_{\fh}}^m$. Moreover,
	we may assume $\widetilde{\xi} \cdot \vac_n$ also lies
	in degrees $\geqslant n \cdot i$.
	With this, the difference $\xi - \widetilde{\xi} \cdot \vac_n$ lies in 
	\[
	F_{i-1}^{PBW}  \bV_{\fh,\kappa_{\fh}}^m \cap
	\bV_{\fh,\kappa_{\fh}}^m
	(\geqslant n\cdot i) \hspace{.5mm} \subset \hspace{.5mm} 
	F_{i-1}^{PBW}  \bV_{\fh,\kappa_{\fh}}^m \cap
	\bV_{\fh,\kappa_{\fh}}^m
	(\geqslant n\cdot (i-1)).
	\]
	\noindent By induction, $\xi - \widetilde{\xi} \in
	\ker(\beta_{n,m})$, and by construction, $\widetilde{\xi} \in
	\ker(\beta_{n,m})$, so we obtain the claim.
	
	Generation is clear as the above calculation for $i = 1$ showed that the 
	first step, i.e. 
	$$F_1^{PBW} \bV_{\fh,\kappa_{\fh}}^m  \cap 
	\bV_{\fh,\kappa_{\fh}}^m(\geqslant n) \hspace{.5mm}= \hspace{.5mm}t^n \fh_O 
	/ t^m \fh_O,$$already 
	generates
	$\ker(\beta_{n,m})$.
	
	\step We now deduce the containment in Lemma \ref{l:alpha-n-m} from the 
	previous step. We remind that $\ft$ denotes the Lie algebra of the Cartan 
	$T$ 
	of $G$.
	
	Recall from \cite{whit} that there
	is a certain level $\kappa_{\ft}$ of $\ft$ such
	that there exist \emph{free-field} homomorphisms
	$$\varphi_n:\sW_{\kappa}^n \to \bV_{\ft,\kappa_{\ft}}^n.$$
	By construction, these are injective morphisms
	that are compatible with loop
	rotation, and by \cite{whit} are strictly
	compatible with the Kazhdan-Kostant and PBW filtrations. Moreover, the 
	free-field morphisms intertwine the transition maps, i.e. one has a 
	commutative diagram
	\[
	\xymatrix{
		\sW_{\kappa}^m \ar@{^(->}[r]^{\varphi_m} \ar[d]_{\alpha_{n,m}}
		& \bV_{\ft,\kappa_{\ft}}^m \ar[d]^{\beta_{n,m}} \\
		\sW_{\kappa}^n \ar@{^(->}[r]^{\varphi_n}
		& \bV_{\ft,\kappa_{\ft}}^n.
	}
	\]	
	Therefore, $\varphi_m(V_{n,m}) \subset W_{n,m} \subset
	\ker(\beta_{n,m})$, which implies that
	$V_{n,m} \subset \ker(\varphi_n\alpha_{n,m}) = \ker(\alpha_{n,m})$, as 
	desired.
	
	\step It remains to
	show generation of $\ker(\alpha_{n,m})$ by $V_{n,m}$.
	
	Let $f \in \fb^-$ denote a principal nilpotent element
	of degree $-1$ with respect to $\check{\rho}$
	and let $\fc$ denote the scheme $(f+\fb)/N$.
	We consider $\bG_m$ acting on $\fc$ in the standard
	way: the action of $\bG_m$ on $\fg$ defined by
	$$\xi \mapsto \lambda \Ad_{\check{\rho}(\lambda)}(\xi), \quad \quad 
	\text{for} \quad \lambda \in \bG_m, \xi \in \fg$$
	induces an action on $\fc = (f+\fb)/N$. It is well-known
	from work of Kostant that one has a $\bG_m$-equivariant isomorphism $$\fc 
	\cong \prod_{i=1}^{\on{rank} \fg } \bA^1,$$ where, if we write $d_i$ 
	for the $i$th exponent of $\fg$, $\bG_m$ acts
	on the $i$th factor of $\mathbb{A}^1$ by the $d_i$th power of the standard
	homothety action.
	
	As $\fc$ is an affine scheme,
	$\fc_O \subset \fc_F$ is a closed subscheme.
	In addition, $\fc_F$ has an action of $(\bG_m)_F$,
	which arises by looping the above action.
	We consider the twist \begin{equation} \label{e:kosslice}\fc_F dt := 
	(\fc/\bG_m)_F 
	\times_{(\mathbb{B}\bG_m)_F}
	\Spec k\end{equation}where $\Spec k \to (\bB \bG_m)_F$ corresponds to
	the line bundle of 1-forms on $\Spec F$. So after
	a choice of non-vanishing 1-form, we may identify $\fc_F dt \simeq \fc_F$,
	but the action of loop rotations is slightly modified.
	We use the similar notation $\fc_O dt$.
	
	Now recall from \cite{whit} that there is a canonical
	isomorphism
	\begin{equation}\label{eq:gr-wn}
	\eta_n:\gr_{\bullet}^{KK} \sW_{\kappa}^n \isom
	\Fun(t^{-n} \cdot \fc_O dt).
	\end{equation}
	\noindent Here $\Fun$ indicates
	the algebra of functions, and
	$t^{-n} \cdot \fc_O dt \subset \fc_F dt$ is
	obtained by acting via $t^{-n} \in (\bG_m)_F(k)$. By construction, this
	isomorphism $\eta_n$ identifies the
	grading on $$\gr_{\bullet}^{KK} \sW_{\kappa}^n = \bigoplus_i \gr_i^{KK}
	\sW_{\kappa}^n$$with the grading on
	$\Fun(t^{-n} \cdot \fc_O dt)$ coming from the action
	of $\bG_m \subset (\bG_m)_O$ on $t^{-n} \fc_O dt$.
	Moreover, $\eta_n$ is equivariant for coordinate changes
	on the disc, and in particular, for loop rotation.
	Finally, for $m \geqslant n$, $\gr_{\bullet}^{KK}(\alpha_{n,m})$
	corresponds to the natural restriction map
	$$\Fun(t^{-m}\cdot \fc_O dt) \onto
	\Fun(t^{-n} \cdot \fc_O dt).$$
	
	From Kostant's description of $\fc$, we
	deduce that $\Ker(\gr_{\bullet}^{KK}(\alpha_{n,m}))$
	is generated by elements in $\gr_i^{KK}(\alpha_{n,m})$
	for some $1 \leqslant i \leqslant d_{\on{rank} \fg}$ that are homogeneous 
	for
	loop rotation, where we can take the loop rotation degrees
	of the generators to lie in the interval $[in,im)$.
	
	We now argue that the above generators may be obtained as symbols of 
	elements of $V_{n,m}$. To see this, consider more generally any $i 
	\geqslant 1$ and an element $$f \in \ker(\gr_i^{KK}(\alpha_{n,m}))$$ which 
	is 
	homogeneous for loop rotation of degree $d \geqslant in$. Lift $f$ to
	$F_i^{KK} \sW_{\kappa}^m$, and then let
	$$\widetilde{f} \in F_i^{KK} \sW_{\kappa}^m$$be the
	degree $d$ component of this lift. Clearly $\widetilde{f}$
	lifts $f$ as well, and $\widetilde{f} \in V_{n,m}$.
	In particular, $\widetilde{f} \in \Ker(\alpha_{n,m})$
	by Step 2 above.
	
	In other words, symbols of elements of $V_{n,m}$ generate
	$\Ker(\gr_{\bullet}^{KK}(\alpha_{n,m}))$. The generation
	immediately follows.\end{proof}

\subsection{Passage to the topological algebra}

Let $\sW_{\kappa}^{as}$ denote the
topological vector space
$$\sW_{\kappa}^{as} := \underset{n}\lim \sW_{\kappa}^n,$$where the topology
is the inverse limit topology. The corresponding pro-object $$``\underset{n} 
\lim" 
\hspace{.5mm} 
\sW_{\kappa}^n \in
\Pro(\sW_{\kappa}\mod^{\heart})$$corepresents the
forgetful functor $\sW_{\kappa}\mod^{\heart} \to \Vect^{\heart}$
by \cite{whit}.
This implies that $\sW_{\kappa}^{as}$ identifies
with endomorphisms of this forgetful functor, so is
canonically a topological algebra.\footnote{More precisely, an $\overset{\rightarrow}{\otimes}$-algebra
in the language of \cite{mys}.}
The left ideals $$\Ker(\sW_{\kappa}^{as} \to
\sW_{\kappa}^n)$$are open and provide a neighborhood basis
of zero. The evident functor identifies
$\sW_{\kappa}\mod^{\heart}$ with the category of
discrete $\sW_{\kappa}^{as}$-modules.\footnote{Although we do not need it, a 
version of this statement for bounded below derived
categories follows from Proposition 3.7.1 of \cite{mys}.}

The Kazhdan--Kostant filtrations and loop rotation gradings extend to the 
topological algebra as follows. For any $i$, define $$F_i^{KK} 
\sW_{\kappa}^{as} := \underset{n} \lim \big(F_i^{KK} 
\sW_{\kappa}^n\big).$$These are 
closed subspaces of 
$\sW_{\kappa}^{as}$ whose union is dense. Similarly, for any $j$ define 
$$\sW_{\kappa}^{as}(\geqslant j) \coloneqq
\underset{n} \lim \big(\sW_{\kappa}^n(\geqslant j)\big).$$

\subsection{Definition of certain ideals}

Let $n \in \bQ^{\geqslant 0}$ be a positive rational number. Define
\[
V_n \coloneqq \sum_{i>0} \big(F_i^{KK}
\sW_{\kappa}^{as} \cap \sW_{\kappa}^{as}
(\geqslant n\cdot i)\big) \subset \sW_{\kappa}^{as}.
\]
\noindent Here we understand
$\sW_{\kappa}^{as}
(\geqslant n\cdot i)$ in the evident way for $n\cdot i$ a rational
number; it coincides with $\sW_{\kappa}^{as}
(\geqslant \ceil{n\cdot i})$. Define $I_n \subset \sW_{\kappa}^{as}$
as the closure of the left ideal generated by
$V_n$.

\begin{pro}\label{p:i_n}
	
	For $n \in \bZ^{\geqslant 0}$, $I_n$ is the kernel of the canonical
	map
	$\alpha_n:\sW_{\kappa}^{as} \to \sW_{\kappa}^n$.
	
\end{pro}

\begin{proof}
	
	Clearly $\ker(\alpha_n) =
	\lim_{m \geqslant n} \ker(\alpha_{n,m})$.
	An element of $V_n$ maps to
	$\lim_{m \geqslant n} V_{n,m}$, so lies in $\ker(\alpha_n)$. Conversely, 
	fix $\xi \in \ker(\alpha_n)$. By Lemma \ref{l:alpha-n-m}, for 
	any $m \geqslant n$ we 
	have $$\xi \in \sW_{\kappa}^{as} \cdot V_n + \ker(\alpha_m).$$
	As the $\ker(\alpha_m)$ form a neighborhood basis of
	zero, this implies that $\xi$ is in the closure
	of the ideal generated by $V_n$.\end{proof}

\begin{cor}
	
	For any $n \in \bQ^{\geqslant 0}$, $I_n$ is open.
	
\end{cor}

\begin{proof}
	
	Clearly $I_m \subset I_n$ for any integer $m \geqslant n$, and
	$I_m$ is open by Proposition \ref{p:i_n}. \end{proof}

\begin{re}\label{r:rat-mods}
	
	We obtain modules $\sW_{\kappa}^n \coloneqq
	\sW_{\kappa}^{as}/I_n$ for all $n \in \bQ^{\geqslant 0}$, generalizing
	the modules from \cite{whit}.
	Each $\sW_{\kappa}^n$ carries an evident
	Kazhdan-Kostant filtrations, and its associated
	graded is the space of
	functions on classical opers with slopes $\leqslant n$.
	
\end{re}

\subsection{}

We will need the following basic results describing
when the ideals $I_n$ ``jump." Below, we write $h$ for the Coxeter number of 
$\fg$, i.e. the maximal exponent $d_{\on{rank} \fg}$. 

\begin{lem}\label{l:in-gen}
	
	For $n \in \bQ^{\geqslant 0}$, define
	\[
	V_n^{\leqslant h} \coloneqq
	\sum_{i=1}^h \big(F_i^{KK}
	\sW_{\kappa}^{as} \cap \sW_{\kappa}^{as}
	(\geqslant n\cdot i)\big) \subset V_n\subset \sW_{\kappa}^{as}.
	\]
	\noindent Then $I_n$ is topologically generated by 
	$V_n^{\leqslant h}$,
	i.e., the closure of the ideal generated by $V_n^{\leqslant h}$ is
	$I_n$.
	
\end{lem}

\begin{proof}
	
	For $n \in \bZ^{\geqslant 0}$, this is immediate from
	the proof of Lemma \ref{l:alpha-n-m} (see
	Step \ref{st:alpha-3}). The same argument applies
	for general $n$, using the general description
	of $\gr_{\bullet} \sW_{\kappa}^n$ from Remark \ref{r:rat-mods}.\end{proof}

\begin{cor}\label{c:ideal-jumps}
	
	Let $h \in \bZ^{>0}$ denote the Coexter number of $G$.
	Then for any $n \in \bQ^{\geqslant 0}$,
	$I_n = I_{\frac{1}{h} \ceil{nh}}$. In particular,
	the ideals $I_n$ are distinct only for
	$n \in \frac{1}{h} \bZ^{\geqslant 0}$.
	
\end{cor}

\begin{proof}
	
	Clearly $V_n^{\leqslant h} = V_{\frac{1}{h} \ceil{nh}}^{\leqslant h}$,
	so we obtain the claim from Lemma \ref{l:in-gen}. \end{proof}

\subsection{Local nilpotence}

In the remainder of this section,
we will show the following result.

\begin{theo}\label{t:whit-nilp}
	
	Let $n$ be a positive integer. Then
	$M \in \sW_{\kappa}\mod^{\heart}$ lies in
	$\Whit^{\leqslant n}(\widehat{\fg}_{\kappa}\mod)^{\heart}$ if and
	only
	if $V_n$ \emph{acts locally nilpotently} on
	$M$, i.e., for every $m \in M$, there exists an
	integer $N \geqslant 0$ such that for every
	$\xi_1,\ldots,\xi_N \in V_n$,
	$\xi_1\ldots\xi_N m = 0$.
	
\end{theo}

Note that Theorem \ref{t:whit-n-sub} is an immediate
consequence of this result. The remainder of this section is devoted to a proof 
of Theorem \ref{t:whit-nilp}. 

\subsubsection{}

We have the following basic observation.

\begin{lem}\label{l:wn-gen}
	
	For an integer $n >0$,
	$\Whit^{\leqslant n}(\widehat{\fg}_{\kappa}\mod)^{\heart}$ is the
	minimal subcategory of $\sW_{\kappa}\mod^{\heart}$ containing
	$\sW_{\kappa}^n$ and closed under colimits and extensions.
	
\end{lem}

\begin{proof}
	
	We have
	$\Whit^{\leqslant n}(\widehat{\fg}_{\kappa}\mod)^{\heart} =
	\widehat{\fg}_{\kappa}\mod^{\o{I}_n,\psi,\heart}$ in the
	notation of \cite{whit}. By construction, the
	module $\sW_{\kappa}^n$ in the left hand side corresponds to 
	$\ind_{\Lie(\o{I}_n)}^{\widehat{\fg}_{\kappa}}(\psi)$.
	Now the claim follows as $\o{I}_n$ is pro-unipotent.\end{proof}

\begin{re}
	
	We can use Remark \ref{r:rat-mods} to give a definition
	of $\Whit^{\leqslant n}(\widehat{\fg}_{\kappa}\mod)$
	and $\Whit^{\leqslant n}(\widehat{\fg}_{\kappa}\mod)^{\heart}$ for
	$n \in \bQ^{\geqslant 0}$: for $n = 0$ it is defined
	as in \cite{whit}, and for $n>0$, we generate under colimits
	using $\sW_{\kappa}^n$. For this definition,
	our argument shows that
	Theorem \ref{t:whit-nilp} is true for $n \in \bQ^{>0}$.
	
\end{re}

\subsection{}

We begin with the following result.

\begin{pro}\label{p:wn-loc-nilp}
	
	For any $n \in \bQ^{>0}$,
	$V_n$ acts locally nilpotently on $\sW_{\kappa}^n$.
	
\end{pro}

\begin{proof}
	
	Fix $v \in \sW_{\kappa}^n$. Choose an integer $i$ such that
	$v \in F_i^{KK} \sW_{\kappa}^n$ and an integer $r$ such that
	$v \in \sW_{\kappa}(\geqslant r)$. We first analyze the action of a single 
	element of $V_n$ on $v$. So, fix 
	$j>0$ and an element $$\xi \in F_j^{KK}
	\sW_{\kappa}^{as} \cap \sW_{\kappa}^{as}
	(\geqslant n\cdot j) \subset V_n.$$
	Clearly $\xi v \in F_{i+j}^{KK} \sW_{\kappa}^n$. In fact,
	we claim that it is in $F_{i+j-1}^{KK} \sW_{\kappa}^n$. To see this, note 
	the	symbol of $\xi$ in $\gr_i \sW_{\kappa}^{as}$
	annihilates the vacuum vector in
	$\gr_{\bullet} \sW_{\kappa}^n$, so as $\gr_{\bullet} \sW_{\kappa}^{as}$
	is commutative and $\gr_{\bullet} \sW_{\kappa}^n$ is generated
	by its vacuum vector, the symbol of $\xi$ annihilates
	all of $\gr_{\bullet} \sW_{\kappa}^n$.
	
	Now take an integer $N$ with $N>i-\frac{r}{n}$ and suppose
	$\xi_1,\ldots,\xi_N \in V_n$.
	We claim that $\xi_1\ldots\xi_N v = 0$. To see this, note from the 
	definition of $V_n$ we may assume that there are
	positive integers $j_k>0$ with each
	$$\xi_k \in F_{j_k}^{KK}
	\sW_{\kappa}^{as} \cap \sW_{\kappa}^{as}
	(\geqslant n\cdot j_k).$$
	Then the above shows that
	$\xi_1\ldots\xi_N v \in F_{i+j_1+\ldots+j_N-N}^{KK} \sW_{\kappa}^n$.
	Clearly this vector is in $\sW_{\kappa}(\geqslant r+n(j_1+\ldots+j_N))$.
	Now observe that
	\[
	r+n(j_1+\ldots+j_N) >
	n(i+j_1+\ldots+j_N-N)
	\]
	\noindent by definition of $N$. The desired vanishing is therefore a 
	consequence of the following Lemma \ref{l:bdded-degs}.	\end{proof}

\begin{lem}\label{l:bdded-degs}
	
	For $n \in \bQ^{>0}$,
	$F_i^{KK} \sW_{\kappa}^n \cap
	\sW_{\kappa}^n (\geqslant ni+1) = 0$.
	
\end{lem}

\begin{proof}
	
	From the description \eqref{eq:gr-wn}
	of $\gr_{\bullet} \sW_{\kappa}^n$,
	we find that $\gr_{j} \sW_{\kappa}^n(m)$ vanishes for
	$m > nj$. Therefore, for $j \leqslant i$,
	$\gr_j \sW_{\kappa}^n(\geqslant ni+1) = 0$, which immediately
	gives the claim.\end{proof}

\subsection{}

We will need one last preparatory lemma. In its statement below, we 
keep the notation of Lemma \ref{l:in-gen}.

\begin{lem}
	\label{l:brackets}
	For $n \in \frac{1}{h}\bZ^{>0}$ and $\xi,\varphi \in V_n$, we have
	$[\xi,\varphi] \in I_{n+\frac{1}{h}}.$
	
\end{lem}

\begin{proof}
	
	We may assume that there are integers $0<i,j \leqslant h$ such that
	$$\xi \in F_i^{KK} \sW_{\kappa}^{as} \cap
	\sW_{\kappa}^{as}(\geqslant ni) \quad \quad \text{and} \quad \quad \varphi 
	\in F_j^{KK} \sW_{\kappa}^{as} \cap
	\sW_{\kappa}^{as}(\geqslant nj).$$
	With this, it follows that $[\xi,\varphi] \in F_{i+j-1}^{KK} 
	\sW_{\kappa}^{as}$,
	as $\gr_{\bullet} \sW_{\kappa}^{as}$ is commutative. Since moreover this 
	element lies in
	$\sW_{\kappa}^{as}(\geqslant n(i+j))$, we have that 
	$$[\xi,\varphi] \in V_{\frac{i+j}{i+j-1} \cdot n}.$$
	As $\frac{i+j}{i+j-1} \cdot n >n$ and $n \in \frac{1}{h}\bZ^{>0}$,
	Corollary \ref{c:ideal-jumps} implies
	$V_{\frac{i+j}{i+j-1} \cdot n} \subset I_{n+\frac{1}{h}}$,
	giving the claim. \end{proof}

\begin{cor}\label{c:nilp}
	
	For $n \in \bQ^{>0}$, let $M$ be a non-zero module
	on which $V_n$ acts locally nilpotently (in the
	sense of Theorem \ref{t:whit-nilp}).
	Then there exists a non-zero
	morphism $\sW_{\kappa}^n \to M$.
	
\end{cor}

\begin{proof}
	
	By Corollary \ref{c:ideal-jumps}, we may suppose $n \in
	\frac{1}{h} \bZ^{>0}$. For any $m \in \frac{1}{h} \bZ^{>0}$, let
	$M_m \subset M$ be the subspace of vectors
	annihilated by $I_m$. By Lemma \ref{l:in-gen},
	this space coincides with the subspace of vectors
	annihilated by $V_m^{\leqslant h}$. Vectors in $M_m$ correspond to maps 
	$\sW_{\kappa}^m \to M$,
	so we equivalently must show that the nonvanishing of $M_n$.

	We will show by descending induction that for any
	$m \in \frac{1}{h}{\bZ}$ with $m \geqslant n$, the subspace
	$M_m \neq 0$. Clearly $$M = \bigcup_m M_m,$$ and so by the nontriviality of 
	$M$ we have the nonvanishing of $M_m$ for $m \gg 0$. Therefore, it suffices 
	to show that if $m \geqslant n$ and $M_{m+\frac{1}{h}} \neq 0$, then 
	moreover $M_m \neq 0$.
	
	We first claim that if $\xi \in V_m$ and
	and $v \in M_{m+\frac{1}{h}}$, then
	$\xi \cdot v \in M_{m+\frac{1}{h}}$. Indeed,
	if $\varphi \in V_{m+\frac{1}{h}} \subset V_m$,
	then $[\varphi,\xi] \in I_{m+\frac{1}{h}}$ by Lemma \ref{l:brackets},
	and so
	\[
	\varphi \xi \cdot v = [\varphi,\xi] \cdot v = 0,
	\]
	as desired. 
	
	It remains to show the nontriviality of $M_m$. Applying Lemma 
	\ref{l:brackets} again, we see that	for $\xi_1,\xi_2 \in V_m$,
	the operators $\xi_i\cdot -:M_{m+\frac{1}{h}} \to M_{m+\frac{1}{h}}$ 
	commute.
	Moreover, as $m \geqslant n$ by assumption, these operators
	are locally nilpotent.
	Therefore, for any finite collection $\xi_1,\ldots,\xi_m \in
	V_n$,
	there exists a non-zero vector in $M_{m+\frac{1}{h}}$
	annihilated by each $\xi_i$. Now the nontriviality of $M_m$ follows from 
	Lemma 
	\ref{l:in-gen}
	and the observation that the inclusion
	$$(V_m^{\leqslant h} \cap I_{m+\frac{1}{h}}) \hspace{.5mm}\subset 
	\hspace{.5mm} V_m^{\leqslant h}$$
	is of finite codimension, which is clear e.g. from
	Lemma \ref{l:bdded-degs}. \end{proof}

\subsection{Conclusion}

In the remainder of the section, we
finish the proof of Theorem \ref{t:whit-nilp}.

\subsection{}

Let $\sA_n \subset \sW_{\kappa}\mod^{\heart}$ be the subcategory
of modules on which $V_n$ acts locally nilpotently. By Proposition 
\ref{p:wn-loc-nilp}, $\sW_{\kappa}^n \in \sA_n$. Clearly $\sA_n$ is closed 
under quotients, extensions, and filtered colimits. Therefore, by Lemma 
\ref{l:wn-gen}, $\Whit^{\leqslant n}(\widehat{\fg}_{\kappa}\mod)^{\heart} 
\subset
\sA_n$.

\subsection{}

Now before proceeding, recall that the
canonical functor $\iota_{n,*}:\Whit^{\leqslant n}(\widehat{\fg}_{\kappa}\mod)
\to \sW_{\kappa}\mod$ is fully faithful (at the derived level),
$t$-exact, and admits a right adjoint $\iota_n^!$. 

We claim that for any $M \in 
\sW_{\kappa}\mod^{\heart}$,
the adjunction map
\[
\varepsilon:
H^0(\iota_{n,*}\iota_n^!(M)) \to M 
\]
\noindent is a monomorphism in $\sW_{\kappa}\mod^{\heart}$. As 
$H^0(\iota_{n,*}\iota_n^!(M)) \in
\Whit^{\leqslant n}(\widehat{\fg}_{\kappa}\mod)^{\heart}$, it
lies in $\sA_n$. Therefore, its subobject
$\Ker(\varepsilon)$ lies in $\sA_n$. We next observe that the map
\begin{equation}
\label{eq:iota-abelian}
\Hom_{\sW_{\kappa}\mod^{\heart}}(\sW_{\kappa}^n,
H^0(\iota_{n,*}\iota_n^!(M))) \to
\Hom_{\sW_{\kappa}\mod^{\heart}}(\sW_{\kappa}^n,
M)
\end{equation}
\noindent is an isomorphism since $\sW_{\kappa}^n \in
\Whit^{\leqslant n}(\widehat{\fg}_{\kappa}\mod)^{\heart}$.
Therefore, $\Hom_{\sW_{\kappa}\mod^{\heart}}(\sW_{\kappa}^n,
\Ker(\varepsilon)) = 0$. By
Corollary \ref{c:nilp}, we find that
$\Ker(\varepsilon) = 0$, as desired.


%
%
%
%

\subsection{}

We now conclude the argument. Suppose $M \in \sA_n$. As above, we need to show 
that the
map
\[
\varepsilon:H^0(\iota_{n,*}\iota_n^!(M)) \to M
\]
\noindent is an isomorphism. We have shown $\varepsilon$ is a monomorphism,
so it remains to see that it is an epimorphism. Again, 
$\Coker(\varepsilon) \in \sA_n$ as it is a quotient of
$M$. The natural map
\[
\Ext_{\sW_{\kappa}\mod^{\heart}}^1(\sW_{\kappa}^n,
H^0(\iota_{n,*}\iota_n^!(M))) \to
\Ext_{\sW_{\kappa}\mod}^1(\sW_{\kappa}^n,
\iota_{n,*}\iota_n^!(M))
\]
\noindent is a monomorphism for cohomological degree reasons.
Also, the map
\[
\Ext_{\sW_{\kappa}\mod}^1(\sW_{\kappa}^n,
\iota_{n,*}\iota_n^!(M)) \to
\Ext_{\sW_{\kappa}\mod}^1(\sW_{\kappa}^n,M)
\]
\noindent is an isomorphism as
$\sW_{\kappa}^n \in \Whit^{\leqslant n}(\widehat{\fg}_{\kappa}\mod)$.
Considering the exact sequence
\[
\Hom_{\sW_{\kappa}\mod^{\heart}}(\sW_{\kappa}^n,\Coker(\varepsilon))
\xar{\partial}
\Ext^1_{\sW_{\kappa}\mod^{\heart}}(\sW_{\kappa}^n,
H^0(\iota_{n,*}\iota_n^!(M))) \to
\Ext^1_{\sW_{\kappa}\mod^{\heart}}(\sW_{\kappa}^n,M)
\]
\noindent we find that the boundary map $\partial$ must
be injective.

Therefore, any map $\sW_{\kappa}^n \to \Coker(\varepsilon)$
lifts to $M$. Applying \eqref{eq:iota-abelian},
we see that any such map is zero. Applying Corollary \ref{c:nilp}
again, we find that $\Coker(\varepsilon) = 0$ as desired.

\section{The baby Whittaker model of monodromic Category $\OO$}
\label{s:moncatowhittaker} \label{s5}
As before, we continue to assume that $\kappa$ is noncritical. In the previous 
section, we studied the adolescent Whittaker filtration on 
$\sW_\kappa\mod^\heartsuit$. Using this, particularly the description of the 
steps in terms of nilpotency of Fourier modes of the generators of 
$\sW_{\kappa}$ obtained in Theorem \ref{t:whit-nilp}, it follows that the 
highest weight representations of $\sW_{\kappa}$ belong 
to the first step in the filtration, i.e. the baby Whittaker category.

In this section, we exactly determine the corresponding objects of the baby 
Whittaker category. That is, we prove an equivalence between a version of 
Category $\OO$ for 
$\sW_\kappa$, whose definition we recall below, and a certain category of 
Iwahori--Whittaker modules for $\gk$, which may not have been considered 
before. 

\subsection{Identification of Category $\OO$}

We begin by reminding the definition of the monodromic category
$\OO$ for $\sW_\kappa$. This was introduced in \cite{lpw}, following
a definition introduced by Arakawa in \cite{ara} in which the
grading was an auxiliary structure, and a definition by Frenkel--Kac--Wakimoto
of postive energy representations \cite{fkw}. Recall that $\sW_\kappa$
contains a canonical conformal vector $\omega$, and that $L_0$
denotes its corresponding energy operator, i.e., the coefficient
of $z^{-2}$ in the field $\omega(z)$.

\begin{defn}
	
	The category $\OO$ is the full subcategory of $\sW_\kappa\mod^\heartsuit$
	consisting of objects $M$ satisfying
	
	\begin{enumerate}
		
		\item Under the action of $L_0$, $M$ decomposes into a sum
		of generalized eigenspaces
		\[
		M = \bigoplus_{d \in \mathbb{C}} M_d, \quad M_d := \{ m \in M:
		(L_0 - d)^N m  = 0, N \gg 0 \}.
		\]
		\noindent
		\item For each $d \in \mathbb{C}$, $M_d$ is finite dimensional
		and $M_{d - n}$ is nonzero for only finitely many $n \in 
		\mathbb{Z}^{\geqslant
			0}$.
	\end{enumerate}
	
\end{defn}

\subsubsection{Identifying a cocomplete variant} \label{sss:icocompletev}Define 
$\OO^{loc}$
to be the full subcategory of $\sW\mod^\heartsuit$ with objects
give by direct limits of objects of $\OO$. As $\OO$ is closed
under subobjects, a module
$M$ belongs to $\OO^{loc}$ if and only if every finitely generated
submodule of $M$ belongs to $\OO$. Therefore,
$\OO^{loc}$ is an abelian category
that is closed under subquotients in $\sW\mod^\heartsuit$.

\begin{re}
	
	If $\kappa$ is negative, then $\OO^{loc}$ may be identified
	with the ind-completion of the full subcategory of finite length
	objects in $\OO$.
	
\end{re}

We now define its counterpart in $\gk\mod^\heartsuit$. Consider
the abelian baby Whittaker category 
\[
{\gk\mod^{\mathring{I}, \psi, \heartsuit}},
\]
\noindent where $\mathring{I} := \mathring{I}_1$. Via restriction, any object 
of this carries 
an
action of $Z$, the center of the universal enveloping algebra
of $\Ad_{t^{-\check{\rho}}} \fg$. Let us consider 
\[
\gk\mod^{\mathring{I}, \psi, Z_f, \heartsuit},
\]
\noindent the full subcategory of $\gk\mod^{\mathring{I}, \psi,
	\heartsuit}$ consisting of objects on which $Z$ acts locally
finitely. I.e., for any vector $v$ in the representation, the
subspace $Zv$ is finite dimensional.

With this, we can state the main result of 
this section. Below, we denote by $$\Psi: \gk\mod \rightarrow \sW_\kappa\mod$$
the functor of Drinfeld--Sokolov reduction.

\begin{theo}
	\label{t:awoo}
	
	The composition $\gk\mod^{\mathring{I},
		\psi, Z_f, \heartsuit} \rightarrow \gk\mod \xar{\Psi[ \langle
		- 2\check{\rho}, \rho \rangle] } \sW_\kappa\mod$ induces an equivalence
	\[\gk\mod^{\mathring{I}, \psi, Z_f, \heartsuit} \simeq \OO^{loc}.\]
\end{theo}

To prove Theorem \ref{t:awoo} we introduce the analogues of Verma modules. 
For a character
$\chi$ of $Z$, consider the $\ag$-module
\[
V_\chi := \mathbb{C}_\chi \underset{Z} \otimes  
\on{ind}_{\Ad_{t^{-\check{\rho}}}
	\fn}^{\ag} \mathbb{C}_\psi.
\]
\noindent We will be interested in its parabolic induction to
$\gk$, namely
\[
\Delta_\chi := \on{pind}_{\ag}^{\gk} V_\chi,
\]
\noindent where as usual `pind' denotes the induction from 
$\Ad_{t^{-\check{\rho}}}
L^+ \fg$ to $\gk$ of the inflation from $\ag$ to $\adt L^+ \fg$
of $V_\chi$.

\begin{lemma}
	\label{l:vzf}
	
	For any $\chi$, the module $\Delta_\chi$
	belongs to $\gk\mod^{\mathring{I}, \psi, Z_f, \heartsuit}.$
\end{lemma}

\begin{proof}
	
	The PBW filtration is a filtration by $\Ad_{t^{-\check{\rho}}}
	\fg$ submodules with associated graded
	\[ \on{Sym}( \gk / \Ad_{t^{-\check{\rho}}} L^+ \fg) \otimes
	V_\chi. \]As the tensor product of an integrable representation
	of $\Ad_{t^{-\check{\rho}}} \fg$ and a $Z$ locally finite
	representation is again $Z$ locally finite, it follows that $Z$ acts 
	locally 
	finitely
	on the associated graded and hence on $M_\chi$.\end{proof}

We also need the following generation result. 

\begin{lemma} \label{l:genbyverms}The  category $\gk\mod^{\mathring{I}, \psi, 
Z_f, \heartsuit}$ 
is the smallest subcategory of  $\gk\mod^{\mathring{I}, \psi, 
\heartsuit}$ which contains the Verma modules $\Delta_\chi$ and is closed under 
colimits and extensions.  
\end{lemma}

\begin{proof} First, note that it 
suffices to show that any object of $\gk\mod^{\mathring{I}, 
\psi, Z_f, \heartsuit}$ may be written as a quotient of a direct sum of finite 
successive extensions of Verma modules. 
To see this, it is enough to show that
	a finitely generated object $N$
	of $\gk\mod^{\mathring{I}, \psi, Z_f, \heartsuit}$ may be written as a 
	quotient of a finite successive extension of Verma modules. Indeed,  Write 
	$\mathring{\mathfrak{i}}$ for the Lie algebra
	of $\mathring{I}$. One has a fiber sequence of Lie algebras
	\begin{equation}\label{e:seslas}
	0 \rightarrow \mathfrak{k}
	\rightarrow \mathring{\mathfrak{i}} \rightarrow 
	\Ad_{t^{-\check{\rho}}}
	\mathfrak{n} \rightarrow 0,
	\end{equation}
	where $\mathfrak{k}$ is $\Ad_{t^{-\check{\rho}}} \mathfrak{k}_1$,
	and $\mathfrak{k}_1$ is the Lie algebra of first congruence subgroup
	of $G_O$. Note that $\Ad_{t^{-\check{\rho}}} \fg$ normalizes
	$\mathfrak{k}$, and $Z$ commutes with $\Ad_{t^{-\check{\rho}}}
	\mathfrak{n}$. Using this, it follows that for any finite dimensional
	subspace $N_0$ of $N$, the image of the action map
	\[U(\mathring{\mathfrak{i}}) \otimes Z \otimes N_0 \rightarrow
	N\]is finite dimensional and stable under the action of $Z$ and
	$\mathfrak{i}$. Moreover, it admits a full flag by subspaces
	which are stable under the action of $Z$ and $\mathring{\mathfrak{i}}$. 
	Taking $N_0$ to contain a set of generators for $N$ and using the 
	parabolic induction of the generated $\Ad_{t^{- \check{\rho}}} \fg$-module 
	yields the claimed surjection.
\end{proof}

We are now ready to prove Theorem \ref{t:awoo}, but first isolate the following 
claim for future reference. 

\begin{lemma} \label{l:v2v} The functor of Theorem \ref{t:awoo} sends Verma 
modules to Verma 
modules, i.e. for any central character $\chi$ we have
	\[
	\Psi[\langle -2 \check{\rho}, \rho\rangle] \hspace{.5mm} \Delta_\chi \simeq 
	M_\chi. 
	\]
\end{lemma}

\begin{proof}This follows from a straightforward
	adaptation from the arguments of Section 4 of \cite{whit}.
	Briefly, the canonical generator of $\Delta_\chi$ defines  
	a vector of highest weight $\chi$ in its
	Drinfeld-Sokolv reduction. This yields
	a map from $M_\chi$. To see this is an isomorphism, note that both
	carry compatible filtrations whose associated gradeds identify
	with the structure sheaf of \begin{equation} 
	\label{e:assgr}\Ad_{t^{-\check{\rho}}}
	\fc_Odt \hspace{.5mm}\subset \hspace{.5mm} \fc_Fdt,\end{equation}
	cf. Equation \eqref{e:kosslice} for the definition of the latter. I.e., if 
	we use the differential of the uniformizer $t$ to identify $\fc_Fdt$ and 
	$\fc_F$, the subscheme in \eqref{e:assgr} consists of loops into $\fc$ that 
	on the factor of $\mathbb{A}^1$ corresponding to the exponent $d_i$ have a 
	pole of order at most $d_i - 1$. 
\end{proof}

\begin{proof}[Proof of Theorem \ref{t:awoo}.] To see that the composition lands 
in $\sW_\kappa\mod^{\heartsuit}$, i.e., in cohomological degree zero, it 
suffices to know the analogous claim for all of $\gk\mod^{\mathring{I}, 
\psi,\heartsuit}$, which was shown in Theorem 7.2.1 of \cite{whit}.
	
	We next show that the image lies in $\OO^{loc}$. We saw in the proof of 
	Lemma \ref{l:genbyverms} that we may characterize 
	\[
\gk\mod^{\mathring{I}, 
	\psi,Z_f, \heartsuit} \subset	\gk\mod^{\mathring{I}, 
		\psi,\heartsuit}
	\]
	as the full subcategory consisting of objects $N$ expressible as a quotient 
	of a direct sum of finite successive extensions of the Verma modules 
	$\Delta_\chi$. By Lemma \ref{l:v2v}, it follows the essential image of 
	$\Psi[\langle - 2 \check{\rho}, \rho \rangle ]$ lies in $\OO^{loc}$.

	It remains to argue that the obtained map $$\Psi: \gk\mod^{\mathring{I},
		\psi, Z_f, \heartsuit} \rightarrow \OO^{loc}$$ is an equivalence.
	As $\gk\mod^{\mathring{I}, \psi, Z_f, \heartsuit} \rightarrow
	\gk\mod^{\mathring{I}, \psi, \heartsuit}$ is a  Serre subcategory
	preserved by taking filtered colimits, it follows from Theorem
	\ref{t:whit-n-sub} that $\Psi$ is fully faithful, and that its essential
	image is closed under taking filtered colimits, quotients, and
	extensions. As every object of $\OO$ admits a filtration, indexed by the 
	nonnegative integers, whose associated graded is a sum of quotients
	of Verma modules, cf. Proposition 6.8 of \cite{lpw}, the essential
	surjectivity follows. \end{proof}

\begin{re}
Having obtained an explicit realization of $\OO^{loc}$ within the baby 
Whittaker category, one may also further ask for a description of the preimage 
of $\OO$. Relatedly, one may ask for a notion of characters in the baby 
Whittaker category which determine the $q$-characters of the associated 
$\sW_\kappa$-modules. We provide answers to both questions in the appendix. 
\end{re}

\subsection{Arakawa's theorem} While we will see a geometric proof of Arakawa's 
results on the `minus' Drinfeld--Sokolov reduction for negative $\kappa$ in
the next section, let us obtain from Theorem \ref{t:awoo} an
algebraic proof here, which applies to any noncritical $\kappa$.

Recall the maximal torus $T$ and Borel $B$ of $G$. Write $B^-$
for the Borel opposite to $B$ with respect to $T$, $I^{-}_0$ for
the corresponding Iwahori subgroup of $G_O$, and $I^-$ for its conjugate
$\Ad_{t^{-\check{\rho}}} I^{-}_0$. Consider the reduction functor
\begin{equation}\label{e:steps}
\Psi_s := \Psi[-2\langle \check{\rho}, \rho \rangle
+ \dim_N]: \gk\mod^{\mathring{I}^-} \rightarrow \sW_\kappa\mod.
\end{equation}
Arakawa showed, following a conjecture of Frenkel--Kac--Wakimoto
\cite{fkw}, that this functor is $t$-exact, sends Verma modules
to Verma modules, and sends simple modules to simple modules
or zero \cite{ara}. As observed in
\cite{whit}, we may factor this through the baby Whittaker category,
i.e.,
\begin{equation}\
\label{e:psis}
\gk\mod^{\mathring{I}^-} \xar{\Av_{\Ad_{t^{-\check{\rho}}}
		N, \psi, *}[\dim_N]}  \gk\mod^{\mathring{I}, \psi} \xar{\Psi[-2
	\langle \check{\rho}, \rho \rangle] } \sW_\kappa\mod.
\end{equation}
As explained in Corollary 7.3.1 of {\em loc. cit.}, the $t$-exactness
is now a general feature of averaging from $N^-$ to $(N, \psi)$
invariants. We now reprove the remaining assertions listed above
via our adolescent Whittaker construction of $\OO$.

Let us introduce relevant notation. For any $\Lambda \in
\ft^*$, where we view $\ft^*$ as the dual of the {\em abstract} 
Cartan,\footnote{The reader may prefer to work with the non-abstract Cartan, 
and should accordingly twist our results by $w_\circ$, the long element of the 
finite Weyl group, due to 
our working with $I^-$.} consider the associated Verma module 
\[ \mathbb{M}_\Lambda :=  \on{ind}_{\Ad_{t^{-\check{\rho}} 
		\fb^-}}^{\Ad_{t^{-\check{\rho}}
		\fg}} \mathbb{C}_\Lambda, \]
\noindent and write $\mathbb{L}_\Lambda$ for its simple quotient. Similarly,
write $\mathbb{M}_{\Lambda, \kappa}$ for the corresponding Verma module
for $\gk$, i.e.,
\[ \mathbb{M}_{\Lambda, \kappa} := \on{pind}_{\Ad_{t^{-\check{\rho}}
		\fg}}^{\gk} \mathbb{M}_\Lambda,\]
\noindent and write $\mathbb{L}_{\lambda, \kappa}$ for its simple quotient.
Finally, write $\lambda$ for the central character of $\mathbb{M}_\Lambda$.
Passing to $\sW_\kappa\mod^\heartsuit$, for $\lambda$ as before
recall that $L_\lambda$ denotes the simple quotient of $M_\lambda$.

\begin{theo}\label{t:arat}
	
	The functor $\Psi_s$, cf. Equation \eqref{e:steps},
	sends Vermas to Vermas, i.e., for any $\Lambda \in \ft^*$, we
	have
	\[ \Psi_s \hspace{.5mm} \mathbb{M}_{\Lambda, \kappa} \simeq M_\lambda. \]
	\noindent Further, $\Psi_s$ sends simples to simples or zero.
	Namely, recalling that $\mathscr{I}$ indexes the finite simple coroots,
	for any $\Lambda \in \ft^*$ we have
	\begin{equation}\label{e:simpszs}
	\Psi_s \hspace{.5mm} \mathbb{L}_{\Lambda, \kappa} \simeq \begin{cases}
	L_\lambda, & \text{if }\langle \Lambda, \halpha_i 
	\rangle \notin
	\mathbb{Z}^{\geqslant 0}, \forall i \in \mathscr{I}, \\ 0, & 
	\text{otherwise.}
	\end{cases}
	\end{equation}
	
\end{theo}

	
	

\begin{proof}
	
	Consider the factorization of $\Psi_s$ through
	the baby Whittaker category as in Equation \ref{e:psis}. We
	saw in Lemma \ref{l:v2v} that the insertion of
	the baby Whittaker category sends $\Delta_\lambda$ to $M_\lambda$,
	and it follows from Theorem \ref{t:awoo} that it exchanges their 
	simple
	quotients. Therefore,
	we are reduced to showing the analogous assertions for the averaging functor
	\[\Av: \gk\mod^{\mathring{I}^-} \xar{\Av_{\Ad_{t^{-\check{\rho}}}
			N, \psi, *}[\dim_N]}  \gk\mod^{\mathring{I}, \psi}. \]
	Let us recall the analog of Theorem \ref{t:arat} in finite type.
	Namely, that
	\[\on{av}: \fg\mod^{N^-} \xar{\Av_{ N, \psi, *}[\dim_N]}  \gk\mod^{N,
		\psi}. \]
	sends $\mathbb{M}_\Lambda$ to the simple object $V_\lambda$, and 
	$\mathbb{L}_\Lambda$
	to $V_\lambda$ or zero, for the exact same conditions on $\Lambda$
	as in Equation \eqref{e:simpszs}. While this is standard, we sketch a proof 
	for the reader's convenience. Indeed, the calculation for Verma
	modules follows from using $t$-exactness and the central character
	of $\mathbb{M}_\Lambda$ to calculate
	\begin{align*}
	H^0 \Hom_{\fg\mod}(V_\lambda, \on{av} \mathbb{M}_\Lambda) & \simeq H^0
	\Hom_{\fg\mod}(\on{ind}_{\fn}^{\fg} \mathbb{C}_\psi, \on{av}
	\mathbb{M}_\Lambda) \\ & \simeq H^{\dim_N} \Hom_{\fn\mod}( \mathbb{C}_\psi,
	\on{Res} \mathbb{M}_\Lambda)\\ & \simeq \det(\fn^*),
	\end{align*}
	\noindent where in the last step one 
	uses that $\mathbb{M}_\Lambda \simeq U(\fn)$
	as $\fn$-modules. The assertion about simples is now immediate
	from $t$-exactness and that the fact that every Verma module
	within a block was sent to the same object.
	
	The assertion about affine Vermas now follows from the equivariance
	of parahoric induction with respect to the categorical action of the Levi, 
	cf. 
	Section A.2 of \cite{tfle}. Namely,
	we have
	\[ \Av \mathbb{M}_{\Lambda, \kappa} = \Av
	\on{pind}_{\Ad_{t^{-\check{\rho}}}\fg
	}^{\gk} \mathbb{M}_{\Lambda} \simeq \on{pind}_{\Ad_{t^{-\check{\rho}}}\fg
	}^{\gk} \on{av} \mathbb{M}_{\Lambda} \simeq 
	\on{pind}_{\Ad_{t^{-\check{\rho}}}\fg
	}^{\gk} V_{\lambda} = M_\lambda.  \]

	For the simples, let us write $\mathbb{A}_{\Lambda, \kappa}$ for the dual
	Verma module of co-highest weight $\Lambda$, i.e., the contragredient
	dual of $\mathbb{M}_{\Lambda, \kappa}$. Similarly, write 
	$\mathbb{A}_\Lambda$ for
	the dual Verma module in $\fg\mod^{B^-}$. For any central character
	$\mu$, we have by adjunction
	\[\Hom_{\gk\mod}( \Delta_{\mu}, \Av \mathbb{A}_{\Lambda, \kappa})  \simeq
	\Hom_{\gk\mod}( \Delta_{\mu}, \mathbb{A}_{\Lambda, \kappa})[\dim_N]\]
	Write $\mathfrak{k}_1$ for the $\check{\rho}$-conjugated first
	congruence algebra, and write the superscript $h \mathfrak{k}_1$
	for Lie algebra cohomology with respect to $\mathfrak{k}_1$.
	Using that $\mathbb{A}_{\Lambda, \kappa}$ is cofree as a $\mathfrak{k}_1$
	module, we may continue
	\[\simeq \Hom_{\fg\mod}( V_\mu , \mathbb{A}_{\Lambda, \kappa}^{h 
		\mathfrak{k}_1})[\dim_N]
	\simeq \Hom_{\fg\mod}( V_\mu, \mathbb{A}_\Lambda)[\dim_N].\]
	From the finite dimensional story, we see this is zero or an
	exterior algebra. In particular, on $\Hom$s in the abelian category
	we have
	\begin{equation}\label{ahom}
	\Hom_{\gk\mod^\heartsuit}(\Delta_\mu, \Av \mathbb{A}_{\Lambda,
		\kappa}) \simeq
	\begin{cases}
	\mathbb{C}, & \text{if }\mu = \lambda, \\
	0, & \text{otherwise}.
	\end{cases}
	\end{equation}

	Consider the factorization of a nonzero map
	\[\Delta_\lambda \onto M \hookrightarrow \Av \mathbb{A}_{\Lambda, 
	\kappa}\]
	through its image. 
	As $M$ is a nonzero quotient of $M_\lambda$, we have a short
	exact sequence
	\[0 \rightarrow K \rightarrow M \rightarrow L_\lambda \rightarrow
	0,\]
	where by mild abuse of notation we write $L_\lambda$ for the simple 
	quotient of $\Delta_\lambda$. We claim $K$ is zero. Indeed, if $K$ is 
	nonzero, then 
	it admits
	a nonzero map $\Delta_\mu \rightarrow K$, for some $\mu \neq \lambda.$
	Composing with the injection $K \rightarrow \Av \mathbb{A}_{\Lambda, 
	\kappa}$,
	we get a contradiction with \eqref{ahom}, as desired. We therefore
	have an image factorization
	\begin{equation}\
	\label{good}
	\Delta_\lambda \twoheadrightarrow L_\lambda
	\hookrightarrow \Av \mathbb{A}_{\Lambda, \kappa}.
	\end{equation}
	If we take the sequence
	\[\mathbb{M}_{\Lambda, \kappa} \twoheadrightarrow L_{\Lambda, \kappa}
	\hookrightarrow \mathbb{A}_{\Lambda, \kappa},\]
	we may apply $\Av$ to obtain a sequence
	\begin{equation} \label{e:gmext}M_\lambda \twoheadrightarrow \Av 
	L_{\Lambda, \kappa} \hookrightarrow
	\Av \mathbb{A}_{\Lambda, \kappa}.\end{equation}
	As we show in the appendix, every object of 
	$\gk\mod^{\mathring{I}, \psi, Z_f, \heartsuit}$
	carries an energy grading, i.e. a locally finite action of (a slight 
	modification of) the Segal--Sugawara energy operator, cf. Section 
	\ref{ss:icato}, and particularly Propositions \ref{p:chv} and \ref{pr111}. 
	Looking at the lowest 
	energy states in \eqref{e:gmext}, we get 
	the analogous sequence
	for $\fg$, i.e., the map is nonzero if and only if $\Lambda$
	satisfies the conditions of \eqref{e:simpszs}. Comparing with \eqref{ahom}
	and \eqref{good} gives the result. \end{proof}

\begin{re} Arakawa also shows that $\Psi_s$ sends dual Verma modules 
to dual Verma modules. However, as 
these may be characterized up to isomorphism by (i) having the same 
Jordan-H\"older content as the corresponding Verma module, and (ii) admitting 
no maps from other Verma modules, this follows from \eqref{ahom}. 
\label{r:costds}
\end{re}

\section{The localization theorem}\label{s6}
Having developed some properties of the adolescent Whittaker filtrations on the 
geometric and representation-theoretic sides of the localization theorem, we 
are ready to turn to its proof.

\subsection{The fully faithful embedding} As discussed in the introduction, the 
existence of a fully faithful embedding follows combining (i) a sufficiently 
robust form of Kashiwara--Tanisaki localization, which incorporates 
renormalized DG categories and categorical loop group actions with (ii) the 
relation between $\sW_{\kappa}$ representations and the Whittaker model of 
Kac--Moody representations established by the second named author in 
\cite{whit}.  

\subsubsection{Kashiwara--Tanisaki localization} Let us state the form of (i) 
which we will need, and along the way fix convenient normalizations for 
discussion of monodromy. 

\subsubsection{}\label{ss:monodromy}Recall the Iwahori subgroup $I^- \subset 
\on{Ad}_{t^{-\check{\rho}}} G_O$. Let us present the affine flag variety as
\[
\on{Fl} = G_F / I^-.
\]

It is well known that the twisting $\kappa$ admits a unique trivialization on 
$\on{Ad}_{t^{-\check{\rho}}} G_O$, and we use this to trivialize its 
restriction to $I^-$. In particular, for an element $\lambda \in \ft^*$ in 
the dual abstract Cartan, we may form 
\[
 D_{\kappa, \lambda-mon}(\Fl) := D_{\kappa}(G_F)^{I^-, -\lambda-mon}.
\]
Explicitly, this is the full subcategory of 
$D_\kappa(G_F)$ generated under shifts and colimits by the essential image of 
the forgetful map from the equivariant category
\[
\on{Oblv}: D_{\kappa}(G_F)^{I^-, -\lambda} \rightarrow D_\kappa(G_F).
\]

\subsubsection{} 

As in Section 7.11 of \cite{raskingl2} (and following \cite{bdh}), we 
\emph{define} a global sections functor
\begin{equation}\label{eq:gamma}
 \Gamma(\Fl, -): D_{\kappa, \lambda-mon}(\Fl) \to  \gk\mod
\end{equation}
\noindent as convolution with the 
{\em monodromic Verma module}
\[
   \mathbb{M}_{\widehat{\lambda}} \in \gk\mod^{I^-, \lambda-mon}. 
\]
\noindent This functor is manifestly 
$D_{\kappa}(G_F)$-equivariant. 

\subsubsection{} We may now state the relevant form 
of localization for Kac--Moody representations. 

\begin{theo}[Kashiwara-Tanisaki localization] \label{t:ktloc} 
If $\lambda$ is regular of level $\kappa$, then 
\[
\Gamma(\Fl, -): D_{\kappa, \lambda-mon}(\Fl) \to  \gk\mod
\]
is fully 
faithful. Its essential image on $\o{I}$-equivariant
objects is the subcategory of $\gk\mod$ 
compactly generated by Verma modules
\[
\mathbb{M}_{\kappa, \theta} \quad \quad \text{for} \quad \theta \in W \cdot 
 \lambda. 
 \]
\noindent 
Moreover, if $\lambda$ is antidominant of level $\kappa$, 
the functor $\Gamma(\Fl,-)$ is $t$-exact. 

\end{theo}

\begin{re}
Using Section A.10 of \cite{raskingl2}, one can show
that $\Gamma(\Fl,-)$ coincides with the usual global
sections functor considered in 
\cite{ktl} and \cite{fgl}. Therefore, for
antidominant $\lambda$, one
can deduce Theorem \ref{t:ktloc} from these previous works. Alternatively, a 
direct proof
is given in the work of J. 
Campbell and the first author on affine Harish--Chandra bimodules 
\cite{ahc}.  \end{re}

\subsubsection{} Having stated Theorem \ref{t:ktloc}, we may obtain from it a 
geometric construction of $\sW$-algebra representations as follows. For any 
twist $\lambda$, let
	\begin{equation}\label{eq:gamma-w}
	  \Gamma(\Fl, -): \msf{Whit}_{\kappa, \lambda-mon}(\Fl) \rightarrow 
	  \sW_{\kappa}\mod 
	\end{equation}
\noindent denote the functor obtained from \eqref{eq:gamma}
by passing to Whittaker invariants and using the
affine Skryabin equivalence
\[
\msf{Whit}(\gk\mod) \simeq \sW_{\kappa}\mod 
\]
of \cite{whit}.

\begin{theo} \label{t:ffloc}

If $\lambda$ is regular of level $\kappa$, this
functor is fully 
faithful. If $\lambda$ is antidominant of level $\kappa$, then 
this functor is $t$-exact. \end{theo}

\begin{proof} 
Fully faithfulness  
follows formally from Theorem \ref{t:ktloc}
using prounipotence of $N_F$. 

For $t$-exactness, note that by construction the insertions 
	\[
	\msf{Whit}^{\leqslant 1}( D_{\kappa, \lambda-mon})(\Fl) \rightarrow 
	\msf{Whit}( D_{\kappa, \lambda-mon}) \quad \quad \text{and} \quad \quad 
	\msf{Whit}^{\leqslant 1}(\gk\mod) \rightarrow \msf{Whit}(\gk\mod)
	\]
	are $t$-exact. We saw in Corollary \ref{c:gendeg1} that the left-hand map 
	is an equivalence. Therefore, it is enough to check the $t$-exactness on 
	baby Whittaker models. Here the assertion clear, as their forgetful 
	functors to the ambient categories are $t$-exact.
\end{proof}

\begin{re}

One can directly show that \eqref{eq:gamma-w} preserves
compact objects for all $\lambda$. Indeed, this follows
for \eqref{eq:gamma} using \cite{raskingl2} 
Lemma 7.12.1. The similar result follows immediately
for invariants for any prounipotent 
compact open subgroup of $G(K)$.
We now obtain the claim for \eqref{eq:gamma-w} using
the baby Whittaker construction and the argument above.

\end{re}

\subsection{Identification of compact generators} Our next task is to produce 
compact generators for the essential image. 

\subsubsection{Geometric preliminaries} Consider the Bruhat decomposition, i.e. 
the stratification of $G_F$ with strata the Schubert cells
\[
  C_w \subset G_F, \quad \quad \text{for} \quad w \in W. 
\]
A convenient indexing for Schubert cells in our situation is 
\[
   C_w := \on{Ad}_{t^{-\check{\rho}}w_\circ} IwI, \quad \quad \text{for} \quad 
   w \in W. 
\] 
In particular, the closure relation between strata is the usual Bruhat order on 
$W$. 

\subsubsection{} It is standard that a cell $C_w$ supports baby Whittaker 
sheaves, i.e. that $D_{\kappa}(C_w)^{\mathring{I}, \psi \times I^-, 
-\lambda}$ is 
nontrivial if and only if $w$ is of minimal length in $W_{\on{f}}w$. In the 
latter 
case, it 
canonically identifies with $\msf{Vect}$, and we write 
\[
       j_{w, !}^\psi, \quad j_{w, !*}^\psi, \quad \text{and} \quad j_{w, 
       *}^\psi, \quad \quad \text{for} \quad w \in W_{\on{f}} \backslash W,
\]
for the corresponding standard, simple and costandard objects of $D_{\kappa, 
\lambda-mon}(\Fl)^\heartsuit$. 

For any $w \in W$, similarly the category associated with a  single stratum 
$$D_{\kappa}(C_w)^{\mathring{I}^- 
\times I^-, -\lambda}$$identifies with 
$\msf{Vect}$, and we 
write 
\[
    j_{w, !}, \quad j_{w, !*}, \quad \text{and} \quad j_{w, *}, \quad \quad 
    \text{for} \quad w \in W,
\]
for the corresponding standard, simple and costandard objects. If we vary the 
twist $\lambda$, we will use the same notations for the similar objects of the 
baby Whittaker and Iwahori-constructible categories.

\subsubsection{} We will use the following standard assertion about averaging 
between 
these two categories. 

\begin{pro}\label{p:avfv} Consider the averaging functor
\[ \Av := \Av_{{I_1}, \psi, *}[\dim_N]:  D_\kappa({I_1}^-
\backslash G_F / I^-, -\lambda) \rightarrow D_\kappa({I_1},
\psi \backslash G_F / I^-, -\lambda). \]
Then for any $w \in W$, if we again denote by $w$ its image in
$W_f \backslash W$, one has isomorphisms:
\begin{equation}\label{e:wavs}
\Av j_{w, !} \simeq j_{w, !}^\psi, \quad \quad
\Av j_{w, *} \simeq j_{w, *}^\psi.
\end{equation}
\end{pro}
\noindent In the case when $\lambda$ and $\kappa$ are both trivial,
this is well known, cf. Section 2 of Arkhipov--Bezrukavnikov
\cite{ab}.

\begin{proof} Consider for $w \in W$ the convolution map 
	\[
           - \overset{I^-}{\star} -: D_{\kappa}(G_F/I^-, -w\lambda) \otimes 
           D_{\kappa}( I^-, w\lambda \backslash G_F / I^-, -\lambda  ) 
           \rightarrow D_\kappa(G_F / I^-, -\lambda). 	
	\]
It was shown in the proof of Theorem 3.7 of \cite{qfle}, which in turn adapted 
an argument of \cite{ly}, that this satisfies
\[
  j_{e, *}^\psi \overset{I^-} \star j_{w, *} \simeq j_{w, *}^\psi  \quad 
  \text{and} \quad  j_{e, *}^\psi \overset{I^-} \star j_{w, !} \simeq j_{w, 
  !}^\psi. 
\]
Similarly one may show that 
\[
  j_{e, *} \overset{I^-} \star j_{w, *} \simeq j_{w, *}  \quad 
\text{and} \quad  j_{e, *} \overset{I^-} \star j_{w, !} \simeq j_{w, 
	!}. 
\]

Using these, by the commutation of $\on{Av}$ and right convolution, we may 
reduce to the case of $w = e$, which in turn is a standard calculation on the 
finite flag variety. 
\end{proof}

\subsubsection{} We are ready to describe a set of compact generators of the 
image. Recall the Harish--Chandra homomorphism 
\[
 \pi: \ft^* \rightarrow W_{\on{f}} \backslash \ft^*. 
\]
\begin{theo} \label{t:compgens}For a twist $\lambda$ which is regular of level 
$\kappa$, the 
essential image of the fully faithful embedding
	\[
	 \Gamma(\Fl, -): \msf{Whit}_{\kappa, \lambda-mon}(\Fl) \rightarrow 
	 \sW_{\kappa}\mod
	\]
	is compactly generated by the Verma modules 
	\[
	 M_\chi, \quad \quad \text{for} \quad \chi \in \pi( W \cdot \lambda). 
	\]
\end{theo}

\begin{proof}
By construction, we have a commutative diagram 
	\[
	\xymatrix{ D_{\kappa, \lambda-mon}(\Fl)^{\mathring{I}^-} \ar[d]_{\on{Av}} 
	\ar[rr]^{\Gamma(\Fl, -)} & & 
	\gk\mod^{\mathring{I}^-} \ar[d]^{\on{Av}} \\ \msf{Whit}^{\leqslant 
	1}_{\kappa, \lambda-mon}(\Fl) \ar[rr]^{\Gamma(\Fl, -)} & &
	\msf{Whit}^{\leqslant 1}(\gk\mod).} 
	\]
By Proposition \ref{p:avfv}, the 
left functor generates under colimits. Therefore,
it suffices to compute the category
generated under colimits by the essential image of the
induced functor:
\[
D_{\kappa, \lambda-mon}(\Fl)^{\mathring{I}^-} \to
\msf{Whit}^{\leqslant 1}(\gk\mod).
\]
\noindent By the diagram and Theorem \ref{t:ktloc},
this is the subcategory
generated by objects 
\[
\on{Av}(\mathbb{M}_{\kappa, \theta}), \quad \quad \text{for} \quad \theta \in W \cdot \lambda.
\] 
By the proof of Theorem \ref{t:arat},  we have
$\on{Av}(\mathbb{M}_{\kappa, \theta}) =
M_{\pi(\theta)}$, 
giving the claim.
\end{proof}

\subsection{The abelian equivalence at negative level}

In the remainder of this section, we specialize to the case of $\lambda$ 
regular and antidominant of level $\kappa$. Recall we have already seen in this 
case that $\Gamma(\Fl, -)$ is $t$-exact. As we will show, it induces an 
equivalence between the abelian category of Whittaker sheaves and a direct sum 
of blocks of Category $\OO$. 

To simplify notation, let us denote the relevant Verma, simple, and dual Verma 
modules for $\sW_\kappa$ by 
\[
   M_w := M_{\pi( w \cdot \lambda)}, \quad L_w := L_{\pi(w \cdot \lambda)}, 
   \quad \text{and} \quad A_w := A_{\pi(w \cdot \lambda)}, \quad \quad 
   \text{for} \quad w \in W_{\on{f}} \backslash W.  
\]
We will also use the integral Weyl group of $\lambda$, which we denote by 
$W_\lambda$; cf. 
Section \ref{ss:regular-weights} for more discussion. With these preparations, 
we may state the remaining theorem of this section. 
\begin{theo} For $\lambda$ antidominant of level $\kappa$, the global sections 
functor induces an equivalence
	\[
	    \Gamma(\Fl, -): \msf{Whit}_{\kappa, \lambda-mon}(\Fl)^\heartsuit \simeq 
	    \bigoplus_{w \in W_{\on{f}} \backslash W / W_\lambda} \OO^{loc}_w
	\]
where $\OO_w^{loc}$ is an indecomposable summand of $\OO^{loc}$ generated under 
extensions and colimits by 
the simple objects
\[
       L_y, \quad \quad \text{for} \quad y \in W_f\backslash W_f w W_\lambda. 
\]
Moreover, this equivalence exchanges standard, simple, and costandard objects, 
i.e. 
\begin{equation}
\label{e:csss}
\Gamma(\Fl, j_{w, !}^\psi) \simeq M_w, \quad \Gamma(\Fl, j_{w, !*}^\psi) \simeq 
L_w, \quad \text{and} \quad \Gamma(\Fl, j_{w, *}^\psi) \simeq A_w, \quad \quad 
\text{for} \quad w \in W_{\on{f}} \backslash W.  
\end{equation}
\label{t:neglev}	
\end{theo}

We now make two remarks before beginning the proof. 

\begin{re} It is straightforward to deduce Theorem \ref{t:neglev} using 
Arakawa's Theorem \ref{t:arat}, by an argument similar to that of Theorem 
\ref{t:compgens}. However, we would like to provide an independent proof 
of the character formulas at negative level. For this reason, below we will 
only use the easy fact that the Drinfeld--Sokolov reduction $\Psi_s$ sends 
Verma 
modules to Verma modules.\end{re}

\begin{re}Passing to `small' categories, an immediate consequence of Theorem 
\ref{t:neglev} is that global sections induces an equivalence
	\[
	\Gamma(\Fl, -): \msf{Whit}_{\kappa, \lambda-mon}(\Fl)^{\heartsuit, f.l.} 
	\simeq \bigoplus_{w \in W_{\on{f}}\backslash W / W_\lambda} \OO_w,
	\]
	where on the left-hand side we mean the finite length (equivalently, 
	compact) objects of the heart and on the right-hand side we mean the 
	corresponding blocks of $\OO$. We remind the reader that while a general 
	block of $\OO$ contains objects of infinite length, this will not occur for 
	the blocks obtained in our localization theorem at negative level (and 
	indeed, this is a consequence of the theorem). 
\end{re}

\begin{proof}[Proof of Theorem \ref{t:neglev}] For ease of notation, we below 
set $\Gamma(-) := \Gamma(\Fl, -).$ A basic property of $\Gamma(
-)$ under our assumption on $\lambda$ is that one has
	\[
	 \Gamma(j_{w, !}) \simeq \mathbb{M}_{\kappa, w \cdot \lambda} \in 
	 \gk\mod^{\mathring{I}^-}, \quad \quad \text{for} \quad w \in W. 
	\]
Therefore, as in the proof of Theorem \ref{t:compgens}, we have
\[
   \Gamma(j_{w, !}^\psi ) \simeq \Gamma(\on{Av} j_{w, !}) \simeq 
   \on{Av} \Gamma(j_{w, !}) \simeq \on{Av} \mathbb{M}_{\kappa, w \cdot 
   \lambda} \simeq M_w,
\]
which establishes the claim of \eqref{e:csss} for standard objects. As a 
consequence, we deduce that the essential image of $\Gamma(-)$ lies in 
$\OO^{loc}$.

   We now turn to the claim for simple objects. We claim it is enough to show 
   that 
   \begin{equation}
       \on{\Hom}_{\sW_{\kappa}\mod^\heartsuit}(M_\chi, \Gamma(j_{w, !*}^\psi)) 
       \simeq 
       0 \quad \quad \text{for} \quad M_\chi \neq M_w.  \label{e:novhoms}
   \end{equation}
   Indeed, the tautological map $j_{w, !} \twoheadrightarrow j_{w, !*}$ yields 
   a map $M_w \twoheadrightarrow \Gamma(j_{w, !*}^\psi).$ In particular, we 
   obtain 
   an exact sequence
   \[
     0 \rightarrow K \rightarrow \Gamma(j_{w, !*}^\psi) \rightarrow L_w 
     \rightarrow 
     0.
   \]
Applying \eqref{e:novhoms}, we deduce that $K$ vanishes, as desired. 

We now prove \eqref{e:novhoms}. From what we have already shown, it is enough 
to argue that, for any $\chi \notin \pi( W \cdot \lambda)$, one has
\[
  \Hom_{\sW_\kappa\mod^\heartsuit}( M_\chi, \Gamma(j_{w, !*}^\psi)) \simeq 0. 
\]
In fact we claim that for such a $\chi$ and any object $N$ of 
$\msf{Whit}_{\kappa, \lambda-mon}(\Fl)$ we have the vanishing 
\begin{equation}\label{e:splitoff}
 \Hom_{\sW_{\kappa}\mod}(M_\chi, \Gamma(N)) \simeq 0. 
\end{equation}
By the compactness of $M_\chi$, it is enough to show the vanishing for $N = 
M_y$, for any $y 
  \in 
  W_{\on{f}}\backslash W.$
To see this, pick a lift $\tilde{\chi}$ of $\chi$, and use
\begin{align*}
   \Hom_{\sW_{\kappa}\mod}(M_\chi, M_{y}) & \simeq \Hom_{\gk\mod^{\mathring{I}, 
   \psi}}( 
   \Delta_{\tilde{\chi}}, \Delta_{y})\\  & \simeq 
   \Hom_{\gk\mod^{\mathring{I}, 
   		\psi}}( 
   \on{Av} \mathbb{M}_{\kappa, \tilde{\chi}}, \on{Av} \mathbb{M}_{\kappa, y 
   \cdot \lambda}) \\ &\simeq   \Hom_{\gk\mod^{I^-, 
   \tilde{\chi}}}(\mathbb{M}_{\kappa, \tilde{\chi}}, \on{Av}^R \on{Av} 
   \mathbb{M}_{\kappa, y \cdot \lambda})
\end{align*}
The vanishing is now clear, as the functor $\on{Av}^R \on{Av}$ is given by 
convolution with an object of 
\[
    D_\kappa(I^-, \tilde{\chi} \backslash G_F / I^-, - y \cdot \lambda) \simeq 
    0,
\]
where the vanishing of the appearing category follows from our assumption on 
$\chi$. Having proven \eqref{e:csss} for standard and simple objects, the 
costandard objects follow by the argument of Remark \ref{r:costds}.

Finally, it remains to show the essential image is a sum of blocks. However, we 
claim this follows from \eqref{e:splitoff}. Indeed, $\OO^{loc}$ lies in the 
subcategory of $\sW_\kappa\mod$ compactly generated by the Verma modules 
$M_\chi$, cf. Proposition \ref{p:appgenvs}, and \eqref{e:splitoff} implies 
this decomposes as the sum of the subcategory generated by $M_w,$ for $w \in 
W$, and the subcategory generated by the remaining Verma modules $$M_\chi, 
\quad \quad \text{for} \quad \chi \notin \pi(W \cdot \lambda).$$ By 
their indecomposability, the simple modules of $L_\chi,$ for $\chi \notin \pi(W 
\cdot \lambda)$, must lie in the latter category, as they admit a nontrivial 
map from $M_\chi$. In particular, for any $w \in W$, we deduce the desired 
vanishing  
\[
 \Ext^1_{\sW_\kappa\mod}(L_w, L_\chi) \simeq \Ext^1_{\sW_{\kappa}\mod}(L_\chi, 
 L_w) \simeq 0.
\] 
 Having shown the essential image is a sum of blocks, it remains to describe 
 its own decomposition into blocks. However, this was determined in Theorem 3.7 
 of \cite{qfle}, completing the proof. \end{proof}

\begin{re} One can show that $\msf{Whit}_{\kappa, \lambda-mon}(\Fl)$ is the 
renormalized derived category of its heart. Namely, its bounded below part 
canonically identifies with the bounded below derived category of its heart, 
and in general the canonical map 
\[
    D^b( \msf{Whit}_{\kappa, \lambda-mon}(\Fl)^{\heartsuit, f.l.}) \rightarrow 
    \msf{Whit}_{\kappa, \lambda-mon}
\]
exhibits the latter as the ind-completion of the former. It follows from this, 
combined with Theorems \ref{t:ffloc} and \ref{t:neglev}, that the subcategory 
of $\sW_{\kappa}\mod$ generated by the Verma modules $M_w, w \in W$, identifies 
with the renormalized derived category of its heart. 
This partially adresses a question
raised in \cite{lpw}.
\end{re}

\appendix

\section{Energy gradings in the baby Whittaker model}

It follows from Theorem \ref{t:awoo} that $\OO$ is equivalent
to subcategory of  $\gk\mod^{\mathring{I}, \psi, Z_f, \heartsuit}$.
In this appendix we will explicitly characterize this subcategory by a positive
energy condition. To account for the conjugation by $t^{-\check{\rho}}$, we as 
usual replace the Segal--Sugawara energy operator with the energy operator
of a different conformal vector, as we explain next.

\subsection{Spectral flow and the Sugawara construction}
Let $\check{\lambda}$ be an integral coweight,
i.e., a cocharacter of the adjoint torus. One has a corresponding
spectral flow automorphism $\Ad_{t^{\check{\lambda}}}$ of
$\gk$. Explicitly, for $X \in \fg$, write $X_n$ for the operator
$X \otimes t^n$ of $\gk$ and consider the corresponding field
$$X(z) = \sum_n X_n z^{-n-1}.$$Using the root decomposition $\fg
= \ft \oplus \bigoplus_{\alpha \in R} \fg_\alpha$, the fields
transform as
\[
\begin{gathered}
\Ad_{t^{\check{\lambda}}} h(z) = h(z) + \frac{\kappa( h,
	\check{\lambda})}{z}, \quad \quad \text{ for } h \in \ft, \\
\Ad_{t^{\check{\lambda}}}  X_\alpha(z) =  z^{\langle \alpha,
	\check{\lambda}\rangle } X_\alpha (z), \quad \quad \text{for } X_\alpha
\in \fg_\alpha.
\end{gathered}
\]
\noindent This induces an automorphism $\Ad_{t^{\check{\lambda}}}$
of the completed enveloping algebra $U_\kappa( \gk)$, which we
denote by the same symbol.

\begin{re} Note that our $\Ad_{t^{\check{\lambda}}}$, which is compatible
with the conventions for the adolescent Whittaker filtration,
differs from that of Arakawa and Frenkel--Kac--Wakimoto by a
sign, i.e., matches their $\Ad_{z^{-\check{\lambda}}}$.\end{re}

\subsubsection{} Consider the usual Segal--Sugawara field $S(z) = \sum S_n 
z^{-n-2}$
of the Kac--Moody vacuum algebra $\mathbb{V}_\kappa$. I.e., writing
$J_a$ for a basis of $\mathfrak{g}$, and $J^a$ for its dual basis
with respect to $\kappa$, we have
\[S(z) =  \frac{\kappa}{2(\kappa - \kappa_c)} \sum_a \normord{J_a(z)J^a(z)}.\]
We presently determine how it transforms under $\Ad_{t^{\check{\lambda}}}$.

\begin{pro}\label{p:dss}
	
	We have the equality of fields
	\[\Ad_{t^{\check{\lambda}}} S(z) = S(z) + \frac{ \check{\lambda}(z)}{z}
	+ \frac{ \frac{1}{2}\kappa(\check{\lambda}, \check{\lambda})}{z^2}.
	\]
	
\end{pro}
\begin{proof}
	
	We must show, for any $n \in \mathbb{Z}$, the equality
	\[\Ad_{t^{\check{\lambda}}} S_n = S_n + \check{\lambda}_n
	+ \delta_{n, 0} \frac{1}{2}\kappa(\check{\lambda}, \check{\lambda}).\]

	We will explain the case of $n = 0$, and the others
	follow by a simpler variant of the same argument. Let $h_i,$
	for $i \in I$, be an orthonormal basis for $\ft$ with respect
	to $\kappa$, and for each positive root $\alpha > 0$ fix vectors
	$e_\alpha$ in $\fg_\alpha$ and $f_\alpha$ in $\fg_{-\alpha}$
	with $\kappa(e_\alpha, f_\alpha) = 1$. Denoting again by $\alpha$
	its image under $\kappa: \ft^* \simeq \ft$, one has $[e_\alpha,
	f_\alpha]  = \alpha$. We have
	\[
	\begin{gathered}
	\Ad_{t^{\check{\lambda}}} S_0= \frac{\kappa}{2(\kappa - \kappa_c)}
	( \sum_{\alpha > 0} \sum_{n \in \mathbb{Z}} \Ad_{t^{\check{\lambda}}}
	\normord{ e_{\alpha, n} f_{\alpha, -n} } + \Ad_{t^{\check{\lambda}}}
	\normord{ f_{\alpha, n} e_{\alpha, -n}} ) \\+ \frac{\kappa}{2(\kappa
		- \kappa_c)}\sum_{i \in I} \sum_{m \in \mathbb{Z}} 
		\Ad_{t^{\check{\lambda}}}
	\normord{h_{i, m} h_{i, -m}}.
	\end{gathered}
	\]
	\noindent Introducing commutators
	to return each sum to its normal order, and recalling that 
	$$\Ad_{t^{\check{\lambda}}}
	h_{i, n} = h_{i, n } + \delta_{n, 0} \kappa( h_i, \check{\lambda}),$$
	we may continue as
	\[
	\begin{gathered}
	\frac{\kappa}{2(\kappa - \kappa_c)}( \sum_{\alpha
		> 0}( 2 \langle \alpha, \check{\lambda} \rangle \alpha_0 + \langle
	\alpha, \check{\lambda} \rangle^2 + \sum_{n \in \mathbb{Z}} \normord{
		e_{\alpha, n} f_{\alpha, -n}} + \normord{ f_{\alpha, n} e_{\alpha,
			-n}})) \\+ \frac{\kappa}{2(\kappa - \kappa_c)}( \sum_{i \in I}
	( 2 \kappa(h_i, \check{\lambda}) h_{i, 0} + \kappa(h_i, \check{\lambda})^2
	+ \sum_{m \in \mathbb{Z}} \normord{h_{i, m} h_{i, -m}})). \\=
	S_0 + \frac{\kappa}{2(\kappa - \kappa_c)}( \sum_{\alpha > 0 }
	2\langle \alpha, \check{\lambda} \rangle \alpha_0 + \langle \alpha,
	\check{\lambda}\rangle^2  + \sum_{i \in I} 2\kappa(h_i, \check{\lambda})
	h_{i, 0} + \kappa(h_i, \check{\lambda})^2 ).
	\end{gathered}
	\]
	\noindent Recalling that $\kappa_c$ is minus one half 
	of the Killing form, we may
	sum over $\alpha$ and $i$ to obtain that the above equals
	\[
	\begin{gathered}
	S_0 + \frac{\kappa}{2(\kappa
		- \kappa_c)}( - \frac{2\kappa_c}{\kappa} \check{\lambda}_0 -
	\kappa_c(\check{\lambda}, \check{\lambda}) + 2 \check{\lambda}_0
	+ \kappa(\check{\lambda}, \check{\lambda}) ) \\ 
	= S_0 + \check{\lambda}_0
	+ \frac{1}{2} \kappa(\check{\lambda}, \check{\lambda}).
	\end{gathered}
	\]

	For the case of $n \neq 0$, one similarly obtains
	\[
	\Ad_{t^{\check{\lambda}}} S_n = S_n + \frac{\kappa}{2(\kappa
		- \kappa_c)}( \sum_{\alpha > 0} 2 \langle \alpha, \check{\lambda}
	\rangle \alpha_n + \sum_{i \in I} 2 \kappa( h_i, \check{\lambda})
	h_{i, n} ) = S_n + \check{\lambda}_n
	\]
	\noindent as desired. \end{proof}

\subsubsection{} Let us use this to determine how $\Ad_{t^{\check{\lambda}}}$
interacts with positive energy conditions on representations.
If $V$ is a representation of $\Ad_{t^{\check{\lambda}}}
\fg $, consider the parabolic induction $\on{pind}_{\Ad_{t^{\check{\lambda}}
	} \fg}^{\gk} V$. It is straightforward to see that the action
of $\Ad_{t^{\check{\lambda}}} S_0$ preserves the image of
the unit $V \rightarrow \on{pind}_{\Ad_{t^{\check{\lambda}}
	} \fg}^{\gk} V$, and acts on it via the element
\[ \frac{\kappa}{2(\kappa - \kappa_c)} \Omega_\kappa,\]
\noindent where $\Omega_\kappa$ is the Casimir operator of 
$\Ad_{t^{\check{\lambda}}
} \fg$ defined with respect to $\kappa$. On the other hand, if
we consider the deformed conformal field $S^{\check{\lambda}}(z)
:= S(z) - \p_z \check{\lambda}(z)$, this again generates a Virasoro
algebra with energy operator
\[S^{\check{\lambda}}_0 = S_0  + \check{\lambda}_0.\]
\noindent  Using Proposition \ref{p:dss} to compare the actions
of $S^{\check{\lambda}}_0$ and $\Ad_{t^{\check{\lambda}}
} S_0$, we obtain the following:

\begin{cor} 
	\label{c:cle}If $V$ is a representation of $\Ad_{t^{\check{\lambda}}
	} \fg$, then $S^{\check{\lambda}}_0$ acts the image of the unit
	$V \rightarrow \on{pind}_{\Ad_{t^{\check{\lambda}} } \fg}^{\gk}
	V$ via the operator of $U(\Ad_{t^{\check{\lambda}} } \fg
	)$ given by
	\[ \frac{\kappa}{2(\kappa - \kappa_c)} \Omega_\kappa - \frac{1}{2}
	\kappa(\check{\lambda}, \check{\lambda}). \]In particular, if
	$\Omega_\kappa$ acts locally finitely on $V$, then $S_0^{\check{\lambda}}$
	acts locally finitely on $\on{pind}_{\Ad_{t^{\check{\lambda}}
		} \fg}^{\gk} V$.
\end{cor}

\subsection{Identifying Category $\OO$}\label{ss:icato}
Specializing
the discussion of the previous subsection to  $\check{\lambda}
= - \check{\rho}$, consider the Kazhdan--Kostant field $S(z)
+ \p_z \check{\rho}(z)$ and write $L_0 := S_0^{-\check{\rho}}$
for its energy operator.

\begin{defn}\label{d:catoiw}
	
	The category $\OO'$ is the full  subcategory
	of $\gk\mod^{\mathring{I}, \psi, \heartsuit}$ consisting of
	objects $M$ satisfying
	
	\begin{enumerate}
		
		\item Under the action of $L_0$, $M$ decomposes into a sum of
		generalized eigenspaces
		\[
		M = \bigoplus_{d \in \mathbb{C}} M_d, \quad M_d := \{ m \in M:
		(L_0 - d)^N m  = 0, N \gg 0 \}.
		\]
		\noindent
		\item For each $d \in \mathbb{C}$, $M_d$ is finite length as
		an $\Ad_{t^{- \check{\rho}}} \fg$ module and $M_{d - n}$
		is nonzero for only finitely many $n \in \mathbb{Z}^{\geqslant
			0}$.
	\end{enumerate}
\end{defn}

For an object $M$ of $\OO'$, if we write $\ell(M_d)$ for
the length of a generalized eigenspace as an $\Ad_{z^{- \check{\rho}}}
\fg$ module, then define its character to be
\[ \ch M := \sum_{d \in \mathbb{C}} \ell(M_d) q^d. \]
In the remainder of this appendix, we first establish some basic properties of 
$\OO'$. We then apply them to prove Theorem \ref{t:awo}, which shows that 
Drinfeld--Sokolov reduction identifies $\OO$ with $\OO'$, and moreover that the 
notions of character on either side are intertwined, up to multiplication by a 
universal 
$q$-series. 

\subsubsection{} We first show that Verma modules belong to our category.
\begin{pro}\label{p:chv}
	
	For any character $\chi$ of $Z$, the
	Verma module $\Delta_\chi$ belongs to $\OO'$ and we have
	\[ \ch \Delta_\chi = q^{ \frac{\kappa}{2(\kappa - \kappa_c)} 
	\Omega_\kappa(\chi)
		- \frac{1}{2} \kappa( \check{\rho}, \check{\rho}) } \prod_{i
		= 1}^\infty \frac{1}{(1 - q^i)^{ \dim \fg} }.\]
\end{pro}
\begin{proof}
	
	As in the proof of Lemma \ref{l:vzf}, consider
	the associated graded of the PBW filtration
	\[ \Sym( \gk / \Ad_{t^{-\check{\rho}}} \fg_O) \otimes
	V_\chi \]
	The energy of the zeroth associated graded is $\frac{\kappa}{2(\kappa
		- \kappa_c)} \Omega_\kappa(\chi) - \frac{1}{2}\kappa( \check{\rho},
	\check{\rho})$ by Corollary \ref{c:cle}. Recalling that $V_\chi$
	is simple, and that for a finite dimensional representation $L$
	of $\Ad_{z^{- \check{\rho}}} \fg$, the tensor product $L
	\otimes V_\chi$ has length $\dim L$, the conclusion follows.\end{proof}

Let us next record the usual relationship between Verma modules
and simple modules, and the usual variant of Jordan-H\"{o}lder
series in the setting of infinite dimensional Lie algebras.

\begin{pro}\label{pr111}
	
	For any character $\chi$ of $Z$, the Verma module
	$\Delta_\chi$ has a unique simple quotient $L_\chi$. Moreover, the
	assignment of $L_\chi$ to $\chi$ yields a bijection between central
	characters and isomorphism classes of simple objects in $\OO'$.
	
\end{pro}

\begin{proof}
	
	For the first assertion, by the simplicity of $V_\chi$,
	it follows that a submodule of $M_\chi$ is proper if and only
	if its intersection with the lowest energy eigenspace is trivial.
	In particular, as every submodule is again $L_0$ locally finite,
	there exists a maximal submodule, and hence a unique simple quotient.
	
	Given an arbitrary simple object $L$, it contains a lowest energy
	subspace, i.e., choose a $d \in \mathbb{C}$ with $L_d$ nonzero
	and $L_{d - n}$ zero for every $n \in \mathbb{Z}^{> 0}$. As $L_d$
	is finite length, we may choose an embedding $V_\chi \rightarrow
	L_d$ of $\Ad_{t^{- \check{\rho}}} \fg$ modules, which induces
	an equivalence $L_\chi \simeq L$. To see that the $L_\chi$, as
	$\chi$ ranges over central characters, are mutually nonisomorphic,
	note that they are distinguished by the action of $\Ad_{z^{-
			\check{\rho}}} \fg$ on their unique lowest energy 
			subspaces.\end{proof}

\subsubsection{} Let us record two further propositions.

\begin{pro}\label{p:22}
	
	Given any scalar $d \in \mathbb{C}$,
	and any object $M$ of $\OO'$, it admits a finite filtration

	\[0 = M^0 \subset M^1 \subset M^2 \subset \cdots \subset M^n
	= M, \]
	where each successive quotient $M^i/M^{i-1}$ is either simple,
	or has no nontrivial subspaces with energy in $d + \mathbb{Z}^{\leqslant
		0}.$
\end{pro}

\begin{pro}\label{p:33}
	
	Any object $M$ of $\OO'$ admits an exhaustive filtration
	\[
	0 = M^0 \subset M^1 \subset M^2 \subset \cdots, \quad \quad
	\underset{n}{\colim} \, M^n \simeq M, 
	\]
	\noindent where for each $n > 0$, the successive
	quotient $M^n/M^{n-1}$ is a direct sum of highest
	weight modules, i.e., quotients of Verma modules.
\end{pro}

As in Proposition \ref{pr111}, one may prove these by adapting
the usual arguments, cf.  Lemma 9.6 of \cite{kitty} and Proposition
6.8 of \cite{lpw}, to our setting. I.e., one considers lowest
energy eigenspaces in lieu of highest weight spaces, and uses
their finite length in lieu of their finite dimensionality.


%
%

\subsubsection{} With these preparations, we are ready to prove the main 
result of this appendix. Recall our notation $\Psi_s$ for the cohomologically 
shifted Drinfeld--Sokolov reduction $\Psi$, namely
\[
\Psi_s = \Psi[ -2\langle\check{\rho}, \rho \rangle].
\]

\begin{theo}\label{t:awo}
	
	The composition $\OO' \rightarrow
	\gk\mod^\heartsuit \xar{\Psi_s} \sW_\kappa\mod$ induces an equivalence
	\[ \OO' \simeq \OO \]
	Moreover, for object $N$ of $\OO'$, its character and that of its 
	Drinfeld--Sokolov reduction are related by
\begin{equation}\
\label{e:ch2w}
\ch \Psi_s N = q^{ \frac{1}{2}\kappa_c(\check{\rho},
	\check{\rho}) + \langle \rho, \check{\rho}  \rangle } \prod_{i
	= 1}^\infty (1 - q^i)^{\dim \fg - \on{rk} \fg} \ch N.
\end{equation}
\end{theo}

\begin{proof} We begin by proving \eqref{e:ch2w}, and start with the Verma 
modules. Recall that
	the character of the Verma module $M_\chi$ for $\sW_\kappa$
	is given by
	\begin{equation} \label{e:chv2}\ch M_\chi = q^{ \frac{\kappa}{2(\kappa - 
	\kappa_c)}
		\Omega_\kappa(\chi) - \frac{1}{2}(\kappa - \kappa_c)(\check{\rho},
		\check{\rho}) + \langle \rho, \check{\rho} \rangle  } \prod_{i
		= 1}^\infty \frac{1}{(1 - q^i)^{\on{rk} \fg}}.\end{equation}
	While this is standard, let us quickly explain how to see the
	energy of its lowest energy eigenspace. Let us calculate $\Psi(\Delta_\chi)$
	using the standard BRST complex. Recall the underlying vector
	space of the latter is $\Delta_\chi \otimes \mathscr{F}$, where 
	$\mathscr{F}$
	is the fermionic ghost vertex algebra associated to $\fn$. If
	we write $\lvert \chi \rangle$ for the canonical generating line
	of $\Delta_\chi$, and $\lvert 0 \rangle$ for the line spanned by the
	vacuum state of $\mathscr{F}$, the lowest energy line of $M_\chi$
	is the image in cohomology of the line
	\[ \lvert \chi \rangle \otimes \det( \Ad_{t^{-\check{\rho}}}
	\fn_O / \fn_O ) \lvert 0 \rangle. \]
	Using Corollary \ref{c:cle} to calculate the energy of the first
	tensor factor and considering the Kazhdan--Kostant grading of
	the determinant, namely $\frac{1}{2}\kappa_c( \check{\rho}, \check{\rho})
	+ \langle \rho, \check{\rho} \rangle$, yields the formula. Comparing 
	\eqref{e:chv2} and Proposition \ref{p:chv} yields \eqref{e:ch2w} for Verma 
	modules. 
	
	 We now deduce
	the general case of \eqref{e:ch2w} as follows. Fix $d \in 
	\mathbb{C}/\mathbb{Z}$
	and consider $\OO'({d})$, the full subcategory of $\OO'$
	consisting of objects whose energy eigenvalues lie in the coset
	$d + \mathbb{Z}$. Similary, consider $\OO(d)$, the full subcategory
	of $\OO$ consisting of objects whose energy eigenvalues lie in
	$$d + \frac{1}{2}\kappa_c(\check{\rho}, \check{\rho}) + \langle
	\rho, \check{\rho} \rangle + \mathbb{Z}.$$By the validity of
	Equation \eqref{e:ch2w} for Verma modules and Proposition \ref{p:33},
	it follows that $\Psi_s$ sends $\OO'(d)$ into $\OO(d)$.
	Moreover, if we pick a coset representative $D \in \mathbb{C}$ of $d$, and 
	consider
	$\OO'(D)$, the full subcategory of $\OO'$ consisting
	of objects whose energy eigenvalues lie in $D + \mathbb{Z}^{\geqslant
		0}$, and similarly define $\OO(D)$, it follows that $\Psi_s$
	sends $\OO'(D)$ into $\OO(D)$. For any nonnegative
	integer $n \in \mathbb{Z}^{\geqslant 0}$, we obtain maps between
	the Serre quotients
	\[ \OO'(D) / \OO'(D + n) \rightarrow
	\OO^{}(D) / \OO^{}(D + n).\]

	By Proposition \ref{p:22}, the Grothendieck group of the left-hand
	category is freely generated as an abelian group by the classes of the 
	Verma modules
	$\Delta_\chi$ with lowest energy lying in $$D, D + 1, \cdots,
	D + n -1.$$Similarly, by Proposition 6.4 of \cite{lpw}
	the Grothendieck group of the right-hand category is freely generated
	by the Verma modules $M_\chi$ with lowest energy in
	$$D_s, D_s+ 1 \cdots, D_s + n -1,$$ where $D_s
	:= D + \frac{1}{2} \kappa_c(\check{\rho}, \check{\rho})
	+ \langle \rho, \check{\rho} \rangle$. Note that the coefficients
	of $$q^{D}, q^{D + 1}, \cdots, q^{D + n - 1}$$
	in the character of any object of $\OO'(D)$ depend only on
	its image in the left-hand Serre quotient, and similarly for
	the coefficients of the first $n$ terms in the $q$-character
	of an object of $\OO(D)$. It therefore follows from the Verma
	case that \eqref{e:ch2w} holds for the first $n$ terms of the
	characters of objects of $\OO'(D)$. Since $n$ and
	$D$ were arbitrary, it follows that \eqref{e:ch2w} holds
	for objects of $\OO'(d)$. As every object of $\OO'$
	is the direct sum of its maximal submodules lying in $\OO'(d)$,
	for $d \in \mathbb{C}/\mathbb{Z}$, the general case of \eqref{e:ch2w}
	follows.
	
	Having shown $\Psi_s$ sends $\OO' \rightarrow \OO$, we deduce
	from Theorem \ref{t:awoo} that it is fully faithful with essential image 
	closed under extensions
	and containing the highest weight modules. To see essential surjectivity,
	we will use an analog of Proposition \ref{p:33} for $\OO$. Namely, for any 
	object $M$ of $\OO$
	we may choose an exhaustive filtration $$M^0 \subset M^1 \subset
	\cdots,$$ whose successive quotients are sums of highest weight
	modules, cf. Proposition 6.8 of \cite{lpw}. As we may assume each summand 
	lies in a fixed term of the associated graded lies in a distinct $\OO(d)$,
	for $d \in \mathbb{C} / \mathbb{Z}$, it follows that each $M^n$
	lies in the essential image of $\OO'$. By \eqref{e:ch2w}, so does the
	direct limit, i.e., $M$,  completing the proof.\end{proof}

\subsubsection{} We finish by proving a statement used in showing the image of 
the localization functor is a sum of blocks at negative level, cf. the proof of 
Theorem \ref{t:neglev}.

\begin{pro} The abelian category $\OO^{loc}$ belongs to the subcategory of 
$\sW_\kappa\mod$ compactly generated by the Verma modules.\label{p:appgenvs}
\end{pro}

\begin{proof} It suffices to show that any object of $\OO'$ belongs to the 
subcategory of $\gk\mod^{\mathring{I}, \psi}$ generated by the Verma modules 
\[
   \Delta_{\chi}, \quad \quad \text{for} \quad \chi \in W_{\on{f}} \backslash 
   \ft^*.
\]
In the notation of the proof of Theorem \ref{t:awo}, we may further reduce to 
the case of an object $N_0$ of $\OO'(D)$. By the proof of Lemma 
\ref{l:genbyverms}, we may find a short exact sequence 
\[
 0 \rightarrow N_1 \rightarrow M_0 \rightarrow N_0 \rightarrow 0
\]
where $M_0$ is a direct sum of finite successive extensions of Verma modules 
which lies in $\OO'(D)$ and $N_1$ is an object of $\OO'(D+1)$. 

Iterating this, we 
obtain a complex 
\[
  \cdots\rightarrow  M_2 \rightarrow M_1 \rightarrow M_0
\]
with each $M_i$ in $\OO'(D+i)$ and of the form above.
	
For any integer $n$, let us write $\sigma^{> n}M$ for its stupid truncation 
$M_{n-1} 
\rightarrow \cdots \rightarrow M_1 \rightarrow M_0.$ Forming the colimit in  
$\gk\mod$ of the stupid truncations under the natural transition maps, it 
suffices to show that the 
obtained map $\colim_n 
\sigma^{> n} M \rightarrow N$
is an equivalence. To check the latter claim, it is enough to see that the 
obtained map 
   \[
       \on{Hom}_{\gk\mod}( \on{ind}_{\mathring{\mathfrak{i}}}^{\gk} 
       \mathbb{C}_\psi, 
       \underset{n}{\colim} \, 
       \sigma^{> n} M) \rightarrow 
       \on{Hom}_{\gk\mod}(\on{ind}_{\mathring{\mathfrak{i}}}^{\gk} 
       \mathbb{C}_\psi, N)
   \]
   is an equivalence. 
   
   To prove this last assertion, let us compute both sides using continuous 
   Chevalley--Eilenberg cochain complexes. The left-hand side produces the 
   direct sum totalization of the bicomplex with underlying vector space 
   \[
       \bigoplus_i \hspace{.5mm} M_i \otimes 
       \on{Sym}(\mathring{\mathfrak{i}}^*[-1]),
   \]
   where the appearing dual is as a topological vector space. We therefore must 
   show that that the natural augmentation to the 
   continuous Chevalley--Eilenberg complex for $N_0$ is an equivalence. To see 
   this note that by 
   construction the energy operator $L_0$ acts locally finitely on the direct 
   sum totalization, and acts on $M_i \otimes 
   \on{Sym}(\mathring{\mathfrak{i}}^*[-1])$ with energy at least $D+i$. Using 
   this, we may safely conclude by applying the standard spectral sequence of a 
   bicomplex, as its convergence is made clear by considering each generalized 
   $L_0$ 
   eigenspace separately. \end{proof}

%
%

%
%
%
%
%
%
%
%
%
%
%
%
%

\bibliographystyle{amsalpha}
\bibliography{samplez}

\end{document}